\documentstyle[11pt,amssymb]{article}
\input{epsf}

\title{ Braided quantum groups 
related to the quantum ``$ax+b$'' group.}

\author{Ma{\l }gorzata Rowicka - Kudlicka\thanks{Supported by KBN grant 
No 2 PO3A 036 18}\\
Institute of Mathematics, Polish Academy of Sciences,\\
\'Sniadeckich 8, 00-950 Warszawa, Poland\\
e-mail:$\;$rowicka@fuw.edu.pl}

\begin{document}

\maketitle
\newcommand{\mqed}{\nopagebreak\centerline{\hfill
\raisebox{3.5ex}[0ex][0ex]{$\Box$}}\\}
\newcommand{\faz}{{\rm Phase\  }}
\newcommand{\ci}{continuous}
\newcommand{\ru}{ unitary representation}
\newcommand{\Qa}{ Quantum 'ax+b' group}
\newcommand{\qa}{ quantum 'ax+b' group}
\newcommand{\MU}{ multiplicative unitary}
\newcommand{\slw}{S.L. Woronowicz}
\newcommand{\tr}{\Delta}
\newcommand{\cH}{{\cal H}}
\newcommand{\K}{{\cal K}}
\newcommand{\cK}{{\cal K}}
\newcommand{\sH}{_{\cal H}}
\newcommand{\sK}{_{\cal K}}
\newcommand{\h}[1]{\hat{#1}}
\newcommand{\lh}{L({\cal H})}
\newcommand{\rf}[1]{{\rm (\ref{#1})}}
\newcommand{\ch}{ {\cal C}(H)}
\newcommand{\chn}{ {\cal C}(H)^N}
\newcommand{\ov}{\overline}
\newcommand{\fh}{F_{\hbar}} 
\newcommand{\hb}{\hbar}
\newcommand{\vt}{V_{\theta}} 
\newcommand{\za}{-\!\!\circ} 
\newcommand{\R}{{\Bbb R}}
\newcommand{\C}{{\Bbb C}}
\newcommand{\Z}{{\Bbb Z}}
\newcommand{\ro}{\rho}
\newcommand{\si}{\sigma}
\newcommand{\be}{\beta}
\newcommand{\de}{\delta}
\newcommand{\ga}{\gamma}
\newcommand{\ta}{\tau}
\newcommand{\N}{{\Bbb N}}
\newcommand{\M}{{\rm M}}
\newcommand{\B}{{\rm B}}
\newcommand{\Lin}{{\rm L}}
\newcommand{\pod}{{\rm d}}
\newcommand{\iz}{\cong}
\newcommand{\eps}{\epsilon}
\newcommand{\fil}{\varphi}
\newcommand{\la}{\lambda}
\newcommand{\qed}{ $\Box$}
\newcommand{\Ci}{C_{\infty}}
\newcommand{\Cir}{C_{\infty}(\R)}
\newcommand{\eh}{e^{\frac{i\hbar}{2}} }
\newcommand{\ehm}{e^{-\frac{i\hbar}{2}} }
\newcommand{\Cgr}{C^{\infty}(\R)}
\newcommand{\Cg}{C^{\infty}}
\newcommand{\Cor}{C_{\rm o}(\R)}
\newcommand{\Co}{C_{\rm o}}
\newcommand{\Cog}{C_{\rm bounded}}
\newcommand{\Cogr}{C_{\rm b}(\R)}
\newcommand{\Lk}{L^{2}}
\newcommand{\Lkr}{L^{2}(\R)}
\newcommand{\po}{\hat{p}}
\newcommand{\qo}{\hat{q}}
\newcommand{\xo}{\hat{x}}
\newcommand{\csta}{$C^{*}$-}
\newcommand{\cstal}{$C^{*}$-algebra \ }
\newcommand{\te}{\otimes}
\newcommand{\ad}{{\rm ad}}
\newcommand{\id}{{\rm id}}
\newcommand{\Mor}{{\rm Mor}}
\newcommand{\Rep}{{\rm Rep}}
\newcommand{\spe}{{\rm Sp }}
\newcommand{\sign}{{\rm sign }\;}
\newcommand{\whe}{\hspace*{5mm}\mbox{\rm where}\hspace{5mm}}
\newcommand{\mand}{\hspace*{5mm} {\rm and} \hspace{5mm}}
\newcommand{\moraz}{\hspace*{5mm} {\rm and} \hspace{5mm}}
\newcommand{\af}{\hspace*{1mm} {\bf \eta} \hspace{1mm}}
\newcommand{\od}{\hspace*{5mm}}
\newcommand{\fu}{{\cal F}}
\newcommand{\fuod}{{\cal F}^{-1}}
\newcommand{\mlot}{\mbox{$\hspace{.5mm}\bigcirc\hspace{-3.7mm}
\raisebox{-.7mm}{$\top$}\hspace{1mm}$}}
\newcommand{\Zak}{\mbox{$-\hspace{-2pt}\comp\,$}}
\newcommand{\dow}{{\bf Proof: }}
%stare
\newcommand{\ut}{\cong}
\newcommand{\Sp}{{\rm sp}}
\newcommand{\infi}{\infty}
\newcommand{\tend}{\rightarrow}
\newcommand{\impl}{\Rightarrow}
\newcommand{\limn}{\lim_{n\rightarrow\infty}}
\newcommand{\limk}{\lim_{k\rightarrow\infty}}
\newcommand{\limt}{\lim_{t\rightarrow\infty}}
\newcommand{\limx}{\lim_{|x|\rightarrow\infty}}
\newcommand{\lime}{\lim_{\eps\rightarrow 0}}
\newcommand{\Lj}{L^{1}}
\newcommand{\Lp}{L^{P}}
\newcommand{\Ls}{L^{S}}
\newcommand{\Li}{L^{\infty}}
\newcommand{\ljr}{{\cal L}^{1}(\R)}
\newcommand{\lk}{ l^{2}}
\newcommand{\lpr}{{\cal L}^{P}(\R)}
\newcommand{\lir}{{\cal L}^{\infty}(\R)}
\newcommand{\Ljr}{L^{1}(\R)}
\newcommand{\Lpr}{L^{P}(\R)}
\newcommand{\Lsr}{L^{S}(\R)}
\newcommand{\Lir}{L^{\infty}(\R)}
\newcommand{\lkn}{ l^2(\N)}
\newcommand{\lp}{{\cal L}^{P}}
\newcommand{\li}{{\cal L}^{\infty}}
\newcommand{\Ljx}{L^{1}(\X)}
\newcommand{\Lkx}{L^{2}(\X)}
\newcommand{\Lpx}{L^{p}(\X)}
\newcommand{\Lqx}{L^{q}(\X)}
\newcommand{\Lix}{L^{\infty}(\X)}
\newcommand{\ljx}{{\cal L}^{1}(\X)}
\newcommand{\lkx}{{\cal L}^{2}(\X)}
\newcommand{\lpx}{{\cal L}^{p}(\X)}
\newcommand{\lix}{{\cal L}^{\infty}(\X)}
\newcommand{\ilkC}{C_{\infty}(\R)\otimes_{C}C_{\infty}(\R)}
\newcommand{\ilk}{\otimes_{\cal C}}
\newcommand{\zw}{CB(\Lkr)}
\newcommand{\zwg}{CB\left(\,L^{2}(G)\,\right)}
\newcommand{\ogr}{B(L^2(\R))}
\newcommand{\hs}{HS(H)}
\newcommand{\spl}{\star}
\newcommand{\milk}{$(i_{1},i_{2},\zw)$}
\newcommand{\mczw}{{\rm Mor(\Cir,\zw)}}
\newcommand{\res}{{\rm Res}}
\newcommand{\pra}{\mbox{${\rm Proj_1}$}}
\newcommand{\prb}{\mbox{${\rm Proj_2}$}}
\newcommand{\Morc}{{\rm Mor}_{\katC}}
\newcommand{\recogr}{{\rm Rep(\Cir,\ogr)}}
\newcommand{\mi}{\hspace*{3mm} {\rm and} \hspace{3mm}}
\newcommand{\mor}{\hspace*{5mm} {\rm or} \hspace{5mm}}
\newcommand{\dla}{\hspace*{5mm} {\rm for} \hspace{5mm}}
\newcommand{\ja}{j_{1}}
\newcommand{\ib}{i_{2}}
\newcommand{\ia}{i_{1}}
\newcommand{\jb}{j_{2}}
\newcommand{\jc}{j_{3}}
\newcommand{\sz}{{\cal S}(\R)}
\newcommand{\gwhom}{$^*$-homomorfizm}
\newcommand{\gw}{$^*$}
\newcommand{\heps}{h_{\eps}} 
\newcommand{\czn}{\Co}
\newcommand{\bh}{B(H) }
\newcommand{\dlad}{for any }
\newcommand{\dlak}{for each }
\newcommand{\Gn}{{\cal G}_n}
\newcommand{\G }{{\rm G}_n}
\newcommand{\rt}{{\bf T}^1}
\newcommand{\jd}{{\bf T}^2}
\newcommand{\Tau}{{\cal T}}

%%%%%

\newcommand{\bfa}{\begin{fakt}}\newcommand{\efa}{\end{fakt}}
\newcommand{\ble}{\begin{lem}}\newcommand{\ele}{\end{lem}}
\newcommand{\bst}{\begin{stw}}\newcommand{\est}{\end{stw}}
\newcommand{\bde}{\begin{defi}}\newcommand{\ede}{\end{defi}}
\newcommand{\bwn}{\begin{wn}}\newcommand{\ewn}{\end{wn}}
\newcommand{\buw}{\begin{uwaga}}\newcommand{\euw}{\end{uwaga}}
\newcommand{\bdy}{\begin{dygresja}}\newcommand{\edy}{\end{dygresja}}
\newcommand{\bwa}{\begin{warning}}\newcommand{\ewa}{\end{warning}}
\newcommand{\bpr}{\begin{przy}}\newcommand{\epr}{\end{przy}}
\newcommand{\btw}{\begin{tw}}\newcommand{\etw}{\end{tw}}
\newcommand{\beq}{\begin{equation}}\newcommand{\eeq}{\end{equation}}
\newcommand{\bit}{\begin{itemize}}\newcommand{\eit}{\end{itemize}}
\newcommand{\bq}{\begin{quote}}\newcommand{\eq}{\end{quote}}
\newcommand{\ba}{\begin{array}}\newcommand{\ea}{\end{array}}

%%%%%%%%%%%%

\newtheorem{defi}{Definition}[section]
\newtheorem{wn}[defi]{Observation}
\newtheorem{tw}[defi]{Theorem}
\newtheorem{lem}[defi]{Lemma}
\newtheorem{aks}{Aksjomat}
\newtheorem{fakt}[defi]{Corollary}
\newtheorem{stw}[defi]{Proposition}
\newtheorem{przy}[defi]{Example}
\newtheorem{uwaga}[defi]{Remark}
\newtheorem{dygresja}[defi]{Dygresja}
\newtheorem{warning}[defi]{Warning}
%********
\begin{abstract} 
We consider quantum group theory on the Hilbert space level. 
We find all unitary representations of three braided quantum groups related 
to the quantum ``$ax+b$'' group. First we introduce an auxiliary 
 braided quantum group, which is  apparently not related to the 
quantum ``$ax+b$'' group, but easy to work with. We find all 
unitary representations of this quantum group. Then we use 
this result to find all unitary representations of another braided 
quantum group, whose $C^*$-algebra is generated by two (out of three) 
generators generating (in the sense of Woronowicz)  
$C^*$-algebra of the quantum 
 ``$ax+b$''- group. We find all unitary representations of this other 
braided quantum group. This is the most difficult 
result needed to classify all unitary representations
 of the quantum ``$ax+b$'' group. 

%This enables us to give a formula 
%for all unitary representations of the quantum ``ax+b'' 
%group in our forthcomming paper \cite{paper2}.

{\bf key words:} unbounded operators -- Hilbert space -- 
 braided quantum group

{\bf MSC-class:} 20G42 (Primary), 47B25 (Secondary).
\end{abstract}

\section{Introduction} 

%In topological quantum group theory a quantum group 
%is  a C$^*$-algebra   equipped with 
% a coassociative comultiplication. 
%From Gelfand-Najmark theorem we know  that every commutative  \cstal is 
% isomorphic  to an algebra of  continuous functions  vanishing at infinity 
%  on some locally compact space. 
%A natural generalisation of what is  to think about 
% non-commutative  \csta algebras as of algebras  of 
%``continuous functions  vanishing at infinity''   
%  on some quantum spaces \cite{wduaC}, \cite{oper}, \cite{connesbook}, 
% where  group operation  is encoded by comultiplication.

%This C$^*$-algebra can be described by  generators and relations \cite{wunb}.
%The generators   
% can be thought of as  non-commutative coordinates on the 'quantum space'.
There are three levels on which one can consider quantum group theory, namely 
the Hopf *-algebra level, the 
\cstal level and the Hilbert space level. In this 
paper we restrict ourselves to  this last level, 
i.e. we consider unbounded operators acting on Hilbert space and encounter various problems related to their domains 
and selfadjoint extensions of  sum of such operators.

 Let $G$ be
a set of closed operators acting on 
 a Hilbert space, invariant under direct sum and unitary tranformations, 
i.e. an operator domain.  
The group structure
may be then introduced by an associative  map from $G\times G$ into 
 $G$. Then, rouhgly speaking, $G$ is quantum group. 
If the Cartesian product $G\times G$ is non-trivial, i.e. operators from  the 
first copy of $G$ do not commute with operators from the second copy of it,
 then  $G$ is a braided quantum group.

% We call a quantum group  compact iff 
% its  algebra   contains identity, in the opposite case we call it  
% a non-compact quantum group. 

Braided quantum groups considered in this paper  are   related 
 to the quantum deformation of the  ``$ax+b$'' group (i.e. group of  affine mappings
of the real line).
They are called  $A$, $N$ and $M$. Operator domain  $A$ is commutative,
 which means 
 that $A$ may be identified with a locally compact space. 
This space is a sum  of three half-lines  with common origin. 
Similarly, $N$ can be identified with a locally compact classical 
space consisting of  four half-lines  with common origin. 
The only quantum, or non-classical,  property of $A$ and  $N$ 
is non-triviality 
 of the Cartesian products  $A\times A$ and $N\times N$. 
The operator domain $M$ is purely quantum -   $M$ is not commutative 
and the Cartesian product $M\times M$ is non-trivial, i.e. there is 
non-trivial braiding.

In  this paper we find all unitary representations 
 of the quantum  groups  $A$, $N$ and $M$. It turns out that every unitary 
representation of  $A$ is  a direct integral of one-dimensional 
  representations.
Since these representations are already found \cite{qexp}, we are 
  able to give a formula for all unitary representation of the braided 
 quantum group  $A$. Having this, it is not difficult to find all
unitary representations of the braided 
 quantum group  $N$. It is done in Section \ref{add}.

The  most interesting result  of this paper is finding all  
unitary representations 
of the braided quantum  group 
$M$. However, the straightforward way to do it requires 
  coping  with operator functions of noncommuting 
unbounded operators.
Functions of unbounded, but commuting, operators are much easier to handle, 
and that is why we start with the case of commutative operator domain $N$.
So it is a great advantage that we invented a simple trick that allows us 
 to ``translate'' our previous
 results concerning $N$ to the case of   $M$.
This way the  problem of  finding  all  unitary representations of 
$M$ is solved in Section \ref{mul}.
This result is essential for classification of  
   all 
unitary representations 
  of the quantum  group "ax+b", which is 
 achieved   in our forthcomming paper 
 \cite{paper2}.

 In the remaining part of this section  we introduce   some non-standard
 notation and  notions used in this paper.
% Szczeg/oln/a uwag/e Czytelnika chcieliby/smy zwr/oci/c na podrozdzia/ly 
%\ref{znalezc}  zawieraj/ace materia/l niestandardowy.

\subsection{Notation}
We  denote Hilbert spaces by   
$ \cH$ and $\cK$, the set of all closed operators acting on $\cH$ by 
 $\ch$, the set of bounded operators by  
 $\B (\cH)$ and the sets of compact and unitary ones 
 by $ CB (\cH)$
 and ${\rm Unit}(\cH)$, respectively.
 The set of all continuous vanishing at infinity functions on a space $X$ will 
be denoted by $\Ci (X)$.
%Przestrzenie Hilberta w tej pracy b/ed/a oznaczane przez $\cH$ i $\cK$. 
We consider only separable Hilbert spaces, usually infinitedimensional.
We  denote scalar product by
$(\cdot|\cdot)$ and it is antilinear in the first variable.
We  consider mainly unbounded linear operators.
All operators considered are densely defined. 
We  use functional calculus of selfadjoint operators 
 \cite{reedI,reedtf,  rudaf}.
We also use the symbol $\sign T$  for partial isometry  
  obtained from  polar decomposition of an selfadjoint  operator $T$.
%We will call an selfadjoint operator  $T$  {\em strictly positive} 
%and write   $T>0$ whenever  $T$ is positive and invertible.
The end of a proof will be marked by  $\Box$.

We  use a non-standard, but very useful notation for ortogonal projections 
and their images  \cite{qexp}, as explained
 below. 
Let $a$ and $b$ be  strongly commuting selfadjoint  
operators  acting on  a Hilbert space  $\cH$. Then by spectral  theorem
there exists a common spectral measure $dE(\lambda)$ such that
\[
a={\int_{\R^2}} \lambda\, dE(\lambda,\mu),
\hspace{1cm}
b={\int_{\R^2}} \mu\, dE(\lambda,\mu)
.\]
For every complex measurable   function    
$f$ of two variables
\[
f(a,b)={\int_{\R^2}} f(\lambda,\lambda')\, dE(\lambda,\lambda').
\]

Let $f$  be a  logical sentence and let   $\chi(f)$
 be 0 if is false, and 1 otherwise.
If ${\cal R}$ is a binary relation on $\R$  
then $f(\lambda,\lambda')=\chi({\cal R}(\lambda,\lambda'))$ is   
 a  characteristic function of a set
$$\Delta=
\{(\lambda,\lambda')\in\R^2:{\cal R}(\lambda,\lambda')\}$$
 and assuming that  $\Delta$ is  measurable $f(a,b)=E(\Delta)$. 
From now on we will write $\chi({\cal R}(a,b))$
 instead of  $f(a,b)$:
\[
\chi({\cal R}(a,b))=
{\int_{\R^2}} \chi({\cal R}(\lambda,\lambda'))\, dE(\lambda,\lambda')
=E(\Delta).
\]
Image of this projector will be denoted  by  ${\cal H}(
{\cal R}(a,b))$, where `$\cH$' is a  Hilbert space, 
 on which operators  $a,b$ act.

Thus  we defined symbols  $\chi(a>b)$, 
 $\chi(a^2+b^2=1)$, \mbox{$\chi(a=1)$}, $\chi(b<0)$, $\chi(a\neq 0)$ 
 etc. They are  ortogonal projections  on appropriate 
  spectral subspaces. For example  ${\cH}(a=1)$ is  
 is an eigenspace of  operator  $a$ for eigenvalue 
 $1$ and   $\chi(a=1)$ is   ortogonal projector on   this  eigenspace.

Generally, whenever  $\Delta$ is a measurable subset of  $\R$,
then  ${\cH}(a\in\Delta)$ is   spectral subspace of an   operator 
$a$ corresponding to   $\Delta$ and  $\chi(a\in\Delta)$ is its   spectral projection.

\subsection{Zakrzewski relation}
 
Let $-\pi<\hb<\pi$.
%Funkcj/e $\sign$ definiujemy nast/epuj/aco:
%\[\sign x=
%\left\{
%\ba{ccc}
%1&\dla&x>0\\
%0&\dla&x=0\\
%-1&\dla&x<0\\
%\ea
%\right.
%\]
The Zakrzewski relation was introduced in  \cite{qexp} by
\bde
Let $R$ and $S$ be selfadjoint operators acting on 
  Hilbert space ${\cH}$. 
Operators  $R$ and $S$ are in    Zakrzewski relation  
 $R\za S$ if
\bit
\item[1.]{$\sign R$ commutes with  $S$ and  $\sign S$ commutes with $R$}
\item[2.]{On subspace $(\ker R)^\perp \cap ( \ker S)^\perp$
we have
\[|R|^{il}|S|^{ik}=e^{i\hb l k} |S|^{ik}|R|^{il}\]
for any $l,k\in\R$. }
\eit
\ede
Observe, that whenever we say that certain $R$ and $S$ satisfy Zakrzewski 
relation we also have to specify parameter $\hb$.

\bpr
Let $\qo$ and  $\po$   denote the position  and momentum operators in 
 Schr\"odinger representation, i.e. we set $\cH=\Lkr$. Then 
 the domain of $\qo$
$$D(\qo)=\{\psi\in \Lkr\;:\; \int_{\R}x^2|\psi(x)|^2dx<\infty\;\}$$
and $\qo$ is multiplication by  coordinate operator on that domain
$$(\qo \psi)(x)=x\psi(x).$$
The domain of  $\po$ consists of all distributions from   $\Lkr$ such that  
$$D(\po)=\{\psi\in \Lkr\;:\; \psi'\in\Lkr\;\} $$
and for any $\psi\in D(\po)$ 
$$(\po f)(x)=\frac{\hbar}{i}\frac{d f(x)}{dx}
,$$
where $-\pi<\hb<\pi$.

Operators  $e^{\po}$ and $e^{\qo}$ acting on $\Lkr$ satisfy  
 Zakrzewski relation. 
\epr

\buw
Two  strictly positive, i.e. positive  and  
 invertible, selfadjoint operators 
 satisfy  Zakrzewski relation 
iff,  they satisfy  Weyl commutation relations \cite{reedtf}.
\euw

\bpr
Let $R$ and $S$ 
 act 
 on Hilbert space 
$\Lkr^{\oplus 4}$. 
%spe/lniaj/ace relacje komutacyjne
%\[
%\label{rs}
%R=R^*,\od S=S^*,\od RS=e^{-i\hb} SR
%\]
% gdzie $-\pi<\hb<\pi$.

%W dalszej cz/e/sci pracy ograniczymy si/e do sytuacji 
%$\hb=\pm \frac{\pi}{2k+3}$, gdzie $k=0,1,2,\dots $.

Let
\[R=\left[\begin{array}{cccccc}
e^{\po} &0&0&0\\0&e^{\po}&0&0\\0&0&-e^{\po}&0
\\0&0&0&\hspace{-2mm}-e^{\po}\\
\end{array}\right],\hspace{10mm}\mand
S=\left[\begin{array}{cccccc}
e^{\qo}&0&0&0\\0&-e^{\qo}&0&0\\0&0&e^{\qo}&0\\0&0&0&-e^{\qo}\\
\end{array}\right],
\]

%Wprowad/zmy oznaczenia ${\cH}(R>0)$ and ${\cH}(R>0)$, gdzie ${\cH}(R>0)$ oznacza podprzestrze/n Hilberta, dla kt/orej $R$ obci/ete do tej przestrzeni is sci/sle dodatnie.. St/ad ${\cH}={\cH}(R>0)\oplus\ker R\oplus  {\cH}(R>0)$.
%Ponadto $S {\cH}(R>0)\subset {\cH}(R>0)$ and $S\ker R\subset \ker R$ and $S {\cH}(R<0)\subset {\cH}(R<0)$
%oraz $R {\cH}(S>0)\subset {\cH}(S>0)$ and $R\ker S\subset \ker S$ and $R {\cH}(S<0)\subset {\cH}(S<0)$.
%Mo/zemy to kr/ocej wyrazi/c pisz/ac, /ze $S$ komutuje z $\sign R$ oraz 
%$R$ komutuje z $\sign S$, 
%gdzie funkcj/e $\sign$ definiujemy nast/epuj/aco:

Hence 
\[\sign R=\left[\begin{array}{cccccc}
1 &0&0&0\\0&1&0&0\\0&0&-1&0\\0&0&0&\hspace{-2mm}-1
\end{array}\right]\hspace{10mm}\mi
\sign S=\left[\begin{array}{cccccc}
1&0&0&0\\0&-1&0&0\\0&0&1&0\\0&0&0&-1
\end{array}\right],
\]
Moreover $|R|$ and $|S|$ 
satisfy  Weyl commutation relations, 
where
\[|R|=\left[\begin{array}{cccccc}
e^{\qo} &0&0&0\\0&e^{\qo}&0&0\\0&0&e^{\qo}&0\\0&0&0&\hspace{-2mm}
e^{\qo}
\end{array}\right]\hspace{10mm}\mi
|S|=\left[\begin{array}{cccccc}
e^{\po}&0&0&0\\0&e^{\po}&0&0\\0&0&e^{\po}&0\\0&0&0&e^{\po}
\end{array}\right]\ .
\]
%Czyli 
%\[R^{il}S^{ik}=e^{ilk\hb}S^{ik}R^{il}\]
%Je/sli $R$ is odwracalne to 
%\[R^{il}S^{ik}R^{-il}=e^{ilk\hb}S^{ik}\]
\epr

All pairs of operators satisfying  Zakrzewski relations are in some sense  
 builded from operators  $e^{\qo}$ and $e^{\po}$.
Precisely, S. L. Woronowicz \cite{qexp} proved 
\bst
\label{vN2}
Let $R$ and $S$  be  operators  
 acting on  Hilbert space 
$\cH$, let  
$\ker R=\ker S=\{0\}$
 and let  $R\za S$.
Then every pair  $(R,S)$   is  
unitarily equivalent to the  
pair $(u\te e^{\po},v\te e^{\qo})$ 
acting on Hilbert space ${\cK}\te \Lkr$, where  
$u,v$ are unitary, 
selfadjoint  and mutually commuting   operators acting on 
 Hilbert space $\cK$
\[ u=u^*=u^{-1},\od v=v^*=v^{-1}\od uv=vu\]
\est

Selfadjoint extensions of sum $R+S$, where $R\za S$, 
proved to be very important in constructing quantum deformation of 
the ``$ax+b$'' group and were studied by Woronowicz  
 in \cite{qexp}. Operator  $R+S$ is symmetric, but in general 
not selfadjoint. To make the thing worse, sometines 
there are even no selfadjoint extensions of such a sum.
However, a selfadjoint extension of the sum $R+S$ exists, 
if there exists a selfadjoint operator $\tau$, such that 
$\tau$ anticommutes with $R$ and $S$ and 
\[ \tau^2=\chi (e^{\frac{i\hb}{2}}RS<0)\;.\] 
Any selfadjoint extension is described uniquely 
by this operator $\tau$,
 so we denote it by $[R+S]_{\tau}$.
 It is given by
\[[R+S]_{\tau}=(R+S)^*|_{D(R+S)+D((R+S)^*)\;\cap\; H(\tau =1)}\;.\]

\subsection{Operator domains and  operator functions}
\label{doifo}

We introduce now two notions very important for understanding this paper,
operator domains and  operator functions.
 They are generalization of  notion of  function  and its domain 
 to the case of function  of  ``non-commuting operator variables''.

The definitions are choosen in such a way that operator functions and 
 domains  respect symmetries of  Hilbert space. 
We are more precise below.

\bde
\label{dop}
Let for any  Hilbert space $\cH$ be given a  subset  
${\cal D}_{\cH}\subset C({\cH})^N$. We say that 
 ${\cal D}$ is 
$N$ - dimensional  {\bf operator domain} if: 
\bit
\item[1.]{For any Hilbert space ${\cH}$ and ${\cK}$ and for any 
 unitary operator $V:{\cH}\rightarrow {\cK}$  
 and  for any   element
$$x=(x_1,x_2,...,x_N)\in {\cal D}_{\cH}\subset \chn$$ we have 
$$V xV^*=(V x_1 V^*,V x_2 V^*, ... ,V x_N V^*)\in {\cal D}_{\cK}$$
}
\item[2.]{
For any space with measure $(\Lambda,\mu)$ and   for any 
  measurable field 
\footnote{Its definition, as well as the definition of a measurable 
field of closed   operators, can be found in  \cite{slwtt}. 
See also \cite{mau}.}
of Hilbert spaces $\{{\cH}(\lambda)\}_{\lambda\in\Lambda}$ and  for any  
 measurable field of closed operators 
$\{{a}(\lambda)\}_{\lambda\in\Lambda}$
 we have 
\[\int_{\Lambda}^\oplus a(\lambda) d\mu (\lambda)
\in D_{ \int_{\Lambda}^\oplus {\cH}(\lambda) d\mu (\lambda)}\]
iff,  $a(\lambda)\in D_{{\cH}(\lambda)}$ for  $\mu-$ almost all  
 $\lambda\in\Lambda$.
}
\eit
\ede

The notion of a measurable field of closed  operators is not widely known, 
but it can be easily reduced to the more popular notion of 
 a measurable field of bounded  operators. 
For any  $T\in\ch$ 
 its {\em $z$-transform} is defined by 
\[z_T=T(I+T^*T)^{-\frac{1}{2}}\ .\]
Observe that  $z_T$ is a bounded operator and $T$  is uniquely 
 determined by  $z_T$. For more details see   \cite{wunb}.

We say that a   field of closed  operators 
\[
\label{pole}
\Lambda\ni\lambda\rightarrow a(\lambda)\in {\cal C}(\cH (\lambda))\]
is measurable  if
\[\Lambda\ni\lambda\rightarrow z_{a(\lambda)}\in {B}(\cH (\lambda))\]
is a measurable field of bounded  operators.
Then there exists a unique  operator
\[a\in {\cal C}(\int_{\Lambda}^\oplus {\cH}(\lambda)d\mu (\lambda))\]
 such that 
\[z_a=\int_{\Lambda}^\oplus z_{a(\lambda)}d\mu (\lambda).\]
We call  operator $a$ a direct integral of the  field \rf{pole} and 
 denote it by  
$\int_{\Lambda}^\oplus {a(\lambda)}d\mu (\lambda).$
For  bounded operators a notion of direct integral introduced above 
 coincides with that used in  \cite{dixmier}.

In particular, when  $\Lambda$ is a countable set  and $\mu$ 
 is a counting measure 
${\cH}=\bigoplus_{\lambda\in\Lambda}{\cH}(\lambda)$
%$, where ${\cH}(\lambda)$ is a Hilbert space
 and condition  2. takes form  
$$\bigoplus_{\lambda\in\Lambda}a(\lambda)\in D_{{\cH}(\lambda)}$$
iff
 $$a(\lambda)\in D_{{\cH}(\lambda)} \od\mbox{\rm for any }\od \lambda\in
\Lambda \ .$$

Observe that   Hilbert space  ${\cH}$ plays a role of a variable in this scheme,
 i.e. difference between an operator domain  $D$ and a set  $D_{\cH}$ is 
 such as between  
 a function $f$ and its value at a point  $x$,  
$f(x)$.
All  operator domains considered   are not closed in the norm topology.

Observe that  an operator domain  is a  category; the  bounded intertwining  
 operators  are its morphisms. For more details see  
 \cite{wduaC} and \cite{wunb}. 
%Theherefore   operator domains considered
 %will be often called  categories.
 
\bpr 
\label{pr:dziedzA}
The set 
\[A_{\cH}=\{ (R,\rho)\in \ch^2\;| \;R=R^*,\; \rho=\rho^*\mi 
R\rho=\rho R\mi \rho^2=\chi (R<0)\}\ \]
is an operator domain. This operator domain was considered in  \cite{qexp}.
\epr

\bpr The set 
\[N_{\cH}=\{ (K,\kappa)\in \ch^2\;| \;K=K^*,\; \kappa=\kappa^*\mi 
K\kappa=\kappa K\mi \kappa^2=\chi(K\neq 0)\}\ \]
is also an operator domain.

% since 
%operatory $K$, $\kappa$ are samosprz/e/zone, to 
% dla ka/zdego unitarnego operatora \mbox{$V:{\cH}\rightarrow K$}  
%operatory 
%$VKV^*$ and $V\kappa V^*$ are r/ownie/z samosprz/e/zone oraz 
% if  $K$ komutuje z $\kappa$, to dla ka/zdego unitarnego operatora 
%$V:{\cH}\rightarrow K$  
%operatory 
%$VKV^*$ and $V\kappa V^*$ b/ed/a komutowa/c. 
%Ponadto suma prosta dw/och samosprz/e/zonych operator/ow is 
%samosprze/zona oraz dla $(K_1,\kappa_1)\in D_{{\cH}_1}$ oraz 
%$(K_2,\kappa_2)\in D_{{\cH}_2}$ sumy proste $(K_1\oplus K_2)$ and $\kappa_1\oplus\kappa_2$ 
%b/ed/a komutowa/c, gdy/z komutuj/a operatory 
%$K_1$ and $\kappa_1$ oraz $K_2$ and $\kappa_2$.
%Podobnie b/edzie dla ca/lek prostych. 
This  operator domain will be discussed  in Section  \ref{add}.
\epr

\bpr 
\label{prM}
The set  
\[M_{\cH}=\{ (b,\beta)\in \ch^2\;| \;b=b^*,\; \beta=\beta^*\mi 
b\beta={\bf -}\beta b
\mi \beta^2=\chi (b\neq 0)\}\ \]
is also an operator domain. This  operator domain will be discussed 
 in Section  \ref{mul}.
\epr

These examples show that  the description of operator domains  is similar 
 to a description of manifold  given by a set of  equations.
For example the sphere 
$$S^2=\{(x_1,x_2,x_3)\in\R^3| x_1^2+x_2^2+x_3^2=1\}\ .$$ 
is described  by giving  {\em coordinates} ($x_1,x_2,x_3$)
and  {\em relations}  between them ($ x_1^2+x_2^2+x_3^2=1$).
Therefore unbounded operators entering descriptions of operator domains
 can be thought of as  ``coordinates on a  quantum space''. 

%${\cal R}$ - relacje komutacyjne elementach z $C({\cH})^N$.
%Relacja to podzbi/or w iloczynie kartezja/nskim element/ow.

The operator functions can be thought of as a  recipee 
what to do with a $N$-tuple of closed 
operators $(a_1, a_2,..., a_N)$ to obtain another closed operator 
 $F(a_1, a_2,..., a_N)$.
\bde
\label{do}
Let $D$ be an  operator domain  and let for any  
Hilbert space  $\cH$ be given a    
 map  $F_{\cH}\;:\;D_{\cH}\rightarrow C({\cH})$.
We say that  $F$ is  a 
 measurable 
{\bf  operator function } if 
\bit
\item[1.]{For any Hilbert spaces ${\cH}$ and ${\cK}$ and for any
 unitary operator $V:{\cH}\rightarrow {\cK}$  
 and for any  element
$$x\in {\cal D}_{\cH}$$ we have 
$$F_{\cK}(V xV^*)=VF_{\cH} (x)V^*$$
}
\item[2.]{
For any  space with a measure $(\Lambda,\mu)$ and
 for any  measurable field of 
 Hilbert spaces $\{{\cH}(\lambda)\}_{\lambda\in\Lambda}$ 
  and for any $a\in D_{\cH}$ having  decomposition 
\[a=\int_{\Lambda}^\oplus a(\lambda) d\mu (\lambda)
\in D_{\int_{\Lambda}^\oplus {\cH}(\lambda) d\mu (\lambda)}\]
 the field of  operators $\{F_{{\cH}(\lambda)}(a(\lambda))\}_{\lambda\in\Lambda}$ 
 is measurable and 
\[F_{\int_{\Lambda}^\oplus {\cH}(\lambda) d\mu (\lambda)}(a)=
\int_{\Lambda}^\oplus F_{{\cH}(\lambda)}(a(\lambda)) d\mu (\lambda).\]}
\eit
\ede

For example, if  $(a_1, a_2, ...,a_N)\in D$ and 
operators $a_i$ are normal, mutually strongly commuting   and 
 for any Hilbert space $\cH$ their joint  spectrum  is contained in 
a set  $\Lambda\subset \C$, then   measurable operator functions on 
  $D$ are simply   measurable  functions on  $\Lambda$. 
It shows that the above definition  is a generalisation of the 
 functional calculus of measurable functions of strongly commuting normal operators to the case of  non-commuting closed, but not necessary normal, 
ones. However in this paper we consider only  operator 
functions of  selfadjoint, but often non-commuting,  operators.

%From the above definition follows  that 
% $F(a_1, a_2,..., a_N)$ commutes with any 
% unitary operator  commuting  with operators $a_1, a_2,..., a_N$.
% Therefore  $F(a_1, a_2,..., a_N)$ is in  bicommutant of  \cstal y generowanej (w sensie wyja/snionym w podrozdziale \ref{not})
% przez operatory $a_1, a_2,..., a_N$. B/edziemy korzystali z tego faktu w rozdziale \ref{sect:ax}.

We also need a notion of an  operator 
maps between   operator domains.

\bde Let $M$  be an  operator domain and let  $N$  be a
$k$ - dimensional operator domain. Moreover, let 
$$F=(F^1, F^2, ..., F^k),$$
 where $F^i$ for  $i=1,2, ...,k$ are operator functions on $M$.    
If for any Hilbert space $\cH$ and for any $m\in M_{\cH}$ 
 we have 
$$F_{\cH}(m)=(F^1_{\cH}(m),F^2_{\cH}(m),...,F^k_{\cH}(m))\in N_{\cH}$$ 
then we call  $F$ an  {\bf operator map} from  an operator  domain 
$M$ into an operator  domain $N$.
\ede 

\subsection{Quantum groups and braided quantum groups on the 
Hilbert space level}

 The definition of  quantum group 
is still under construction  \cite{vaeskrotki}, there are however many 
examples and one knows approximately what a quantum group should be.
%Roughly speaking, 

 Let $G$ be an operator domain and let $G\times G$ denote an operator 
domain
\[G\times G:=\{ (x,y)\;|\; x,y\in G \mi xy=yx\}\]
Let $\cdot$ be an operator map
$$\cdot{ } \;\;\;:G\times G\ni (x,y)\rightarrow xy\in G$$.

Roughly speaking, a quantum  group $G$  is such an operator 
domain $G$ equipped with an associative operator map $\cdot$.

\bpr[Quantum  $SU_{{q}}(2)$ group]
%Kwantowa grupa $SU_q(2)$ 
Let us define an operator domain
$$SU_{{q}}(2)=\left\{
(\alpha,\gamma)\in \bh
:\; 
\ba{c}
\alpha\alpha^*+\gamma\gamma^*=I;\;\\
\alpha\gamma={q}\gamma\alpha;\;\\
\gamma\gamma^*=\gamma^*\gamma;\;\\
\alpha\gamma^*={q}\gamma^*\alpha;\;\\
\alpha^*\alpha+{q}^2\gamma^*\gamma=I\;
\ea\;
\right\}\;.$$

Moreover, let us define an operator map $\cdot$ by
% $SU_{{q}}(2)$ 
\[
\cdot:
SU_{{q}}(2)
\times SU_{{q}}(2)\ni 
((\alpha_1,\gamma_1),(\alpha_2,\gamma_2))
\tend (\alpha,\gamma)\in SU_{{q}}(2)\]
where 
\[\alpha =\alpha_1\alpha_2-_{{q}}{\gamma_1}^*\gamma_2\mi 
\gamma=\gamma_1\alpha_2+{\alpha_1}^*\gamma_2\]

This operation is associative, and $SU_{{q}}(2)$ together with this operation 
 forms a quantum group. Since in the definition of the operation 
domain $SU_{{q}}(2)$ we restricted ourselves to the bounded operators, 
$SU_{{q}}(2)$ is a compact quantum group.
\epr
%************************
%\bpr[Quantum  ``$a{z}+b$'' at roots of unity] 
%$$G_{\rm H}=\left\{\;(a,b,
%{\beta})\in{\cal C}({\rm H})^3 \;
%:\;
%\ba{c} \; \\
% a \mbox{\rm invertible} \\
%a b=qb a
%\ea
%\right\}
% $$

%Group operation in  $G_{\rm H} $
%\[
%%\cdot:
%G_{\rm H}\times 
%G_{\rm H}\ni ((a_1,b_1),
%(a_2,b_2))
%\tend (a,b)\in G_{\rm H}\]
%gdzie 
%\[a =a_1\te a_2\mi 
%b=a_1\te  b_2{\dot{+}} b_1\te I\]
%\epr

%************************
\bpr[Quantum ``$a{x}+b$'' group] 
Let us define an operator domain 
$$G_{\rm H}=\left\{\;(a,b,
{\beta})\in{\cal C}({\rm H})^3 \;
:\;
\ba{c} \; \\
 a>0 \\
a\za b\\
(b,\beta)\in M_{\rm H} \\
a\beta=\beta a
\ea
\right\}\;,
 $$
where $M_{\rm \cH}$ is the operator domain defined in Example \ref{prM}.

Group operation $\mlot$ on  $G_{\rm H} $ is given by
\[
\mlot:
G_{\rm H}\times G_{\rm H} 
\ni ((a_1,b_1,\beta_1),
(a_2,b_2,\beta_2))
\tend (a,b,\beta)\in G_{\rm H}\]
where 
\[a =a_1\te a_2\mi 
b=[a_1\te  b_2+ b_1\te I]_{(-1)^{{k}}(\beta_1\te\beta_2)
\chi(b_1\te b_2<0)}\]
Formula for  $\beta$  is much more complicated.
This additional ``generator'' $\beta$ is needed to ensure the 
existence of a selfadjoint extension of $a_1\te  b_2+ b_1\te I$.

To make group operation $\mlot$ associative one has assume that 
\[ \hb=\pm \frac{\pi}{2{k}+3}, \whe k=0,1,2,\dots
\]
and $k$ is the same as chosen in formula for selfadjoint 
 extension of $a_1\te  b_2+ b_1\te I$.

An operator domain $G$ with operation $\mlot$ defined as above 
forms the quantum ``$ax+b$''  group \cite{ax+b}. 
Observe that from  \mbox{\rm $a\za b$} follows that operators 
$a$ and $b$ are 
 not bounded.
Therefore the quantum ``$ax+b$''  group is non-compact.
\epr
The theory of non-compact quantum  groups is more difficult, more interesting and 
 less developed than that of  
 compact ones.
The most important examples of  non-compact quantum  groups are 
  the quantum $E(2)$ group, quantum Lorentz group  and quantum  
groups "$ax+b$" and  "$az+b$" 
\cite{oe2}, \cite{e2d}, \cite{podleslorenz}, \cite{pusz},
\cite{ax+b}, \cite{az+b}.

What we study in this paper are braided quantum groups.
 The  main difference between a braided quantum group and 
 a quantum group is that a group operation on a braided 
 quantum group $G$ is defined on an operator domain 
\[G^2:=\{ (x,y)\;|\; x,y\in G \mi x,y \
\mbox{\rm satisfy certain relations}\}\;.\]
Usually we do not assume that operators from both copies 
 of $G$ commute, so in general a braided quantum group is 
not a quantum group.
A group operation on a braided quantum group $G$ shoul be the  
operator map 
\[\mlot:G^2\ni (x,y)\rightarrow x\mlot y\in G\;\] 
 which is associative.

Roughly speaking, a braided quantum group $G$ is an operator domain $G$ 
 equipped  with such an  operator map $\mlot$.

The objects $N$, $A$ and $M$ investigated in this paper are 
examples of braided quantum groups.
%********
\section{The braided quantum group  $N$}
\label{add}
Let ${\cal H}$ be a separable, infinitedimensional 
Hilbert space.
Consider  operators  $R $ 
 and   $\rho$  acting on $\cH$ and such that 
\[R=R^*\mand \rho=\rho^*\mand \rho R=R\rho\mand \rho^2=\chi(R\neq 0)\ .  \]
Let
\beq
\label{rkN}
N_{\cal H}=\{ (R,\rho)\in \ch\;:\; R=R^*\;\;\; 
\rho=\rho^*\;\;\rho R=R\rho\;\; \rho^2=\chi(R\neq 0)\;\}.
\eeq
%\bpr Pary z $A_{\cal H}$ to: $(3I,0),(-\qo^2-2I,I),(-I,I),
%(e^{\qo},0),(-e^{\po},I)$.
%\epr
%
It is easily seen that  $N$ is an operator domain.
To ensure existence of 
 selfadjoint extensions of a 
sum $R+S$, where $(R,\rho),(S,\sigma)\in N_{\cH}$, we introduce an additional 
condition for pairs  
$(R,\rho),(S,\sigma)\in N_{\cH}$ to fulfill. This condition 
%we call a  concordance   condition.
% The concordance condition for  pairs $(R,\rho),(S,\sigma)\in N_{\cH}$ 
is
\beq
\label{warzgodN}
R\za S \mi  S\rho=-\rho S \mi 
R\sigma=-\sigma R\mi \rho\sigma=\sigma\rho\ .\eeq
If pairs $(R,\rho),(S,\sigma)\in N_{\cH}$ satisfy 
  this  condition we write 
$$(R,\rho),(S,\sigma)\in N_{\cH}^2\;,$$
i.e.
\[N_{\cal H}^2=\{ (R,\rho),(S,\sigma)\in N_{\cal H}\;:\;
R\za S \mi  S\rho=-\rho S \mi R\sigma=-\sigma R\mi \rho\sigma=\sigma\rho
\;\}\ .\]
%\[N_{\cal H}^2=\{ (R,\rho,S,\sigma)\; :\; (R,%\rho),(S,\sigma)\in N_{\cal H}
%\mi R\za S \mi R\sigma =-\sigma R\mi 
%S\rho=-\rho S\mi \rho\sigma =\sigma\rho \}\]
Observe, that  $N^2$ is also an operator domain.
 
To give a formula for selfadjoint extensions of a 
sum $R+S$ we need the quantum exponential function.

\subsection{The special function  $\vt$ and 
 the quantum exponential function $\fh$}
 The special function  $\vt$ is defined by 
\[\vt(x)=\exp\left\{\frac{1}{2\pi i}\int_0^\infty \log(1+a^{-\theta})
\frac{da}{a+e^{-x}}\right\}\]
 for any  $x\in \C$ such that  $|\Im x|<\pi$.
$\vt$ can be extended to a function meromorphic on  $\C$.
Let \[
\Omega^{+}_{\hbar}=\left\{\;r\in\C\setminus \{0\}\;:\;\arg r \in [0,\hbar]\;
\right\}
\]
 \[
\Omega^{-}_{\hbar}=\left\{\;r\in\C\setminus \{0\}\;:\;\arg r \in [-\pi,
\hbar-\pi]\;
\right\}
\]
\[\Delta=\Omega^{+}_{\hbar}\times \{0\}\cup\Omega^{-}_{\hbar}\times
\{-1,1\} \]
The quantum exponential function   $\fh$ is a  $q$-analogue of the 
  exponential function fit 
 for  operators satisfying commutation rules of the type  \rf{warzgodN}.
 What we mean  exactly by this is explained in Proposition 
\ref{st:oexp}.
The  function  $\fh$ is defined 
  for any $(r,\rho)\in \Delta$ by
\[
\fh(r,\rho)=[1+i\rho (-r)^{\frac{\pi}{\hbar}}]\vt(\log r)
\]
In particular, for  $(r,\rho)\in \Delta_{real}:=\{(r,\rho)\in \Delta:r\in \R\}$
\beq
\fh(r,\rho)=\left\{
\begin{array}{ccc}
\vt (\log r)& \mbox{\rm dla } & r>0 \;\;\;\mbox{\rm  and } \rho = 0\\
\{1+i\rho |r|^{\frac{\pi}{\hbar}}\}\vt(\log |r|-\pi i)& 
\mbox{\rm dla} & r<0 \;\;\;\mbox{\rm and } \rho = \pm 1
\end{array}
\right.
\eeq
Axiomatic introduction of $\fh$ and some other properties of  
$\fh$ and $\vt$ can be found in Section 2 of \cite{qexp}. 
The function  $\fh$ is extended to the closure of  $\Delta_{real}$ by setting 
\[\fh (0,\rho)=1.\]
\subsection{The group operation on $N$}
Observe that  if  $(R,\rho)\in N_{\cal H}$ then 
 $R$ commutes with  $\rho\chi(R<0)$ and the  joint spectrum
 of these operators is the closure of   $\Delta_{real}$.
Therefore $\fh(R,\rho\chi(R<0))$ is well defined.
Moreover, since $(R,\rho\chi(R<0)) A_{\cal H}$ and 
$$\fh:\Delta_{real}\rightarrow \C$$ is a measurable function, 
we know from Introduction  that
 $\fh$ is an operator function defined on $A$.

 Let  $\left( (R,\rho),(S,\sigma)\right)
\in N_{\cal H}^2 $. 
Assume that $\ker S=\{0\}$. This asumption is not very restrictive 
since every  $S$ is a direct sum of 
  invertible  $S_1$   and $S_2=0$, and the case  $S_2=0$ 
 is trivial. 
Define 
\beq
\label{ttau}
T=e^{\frac{i\hbar}{2}}S^{-1}R\;\;\;\tau=(-1)^k\rho\sigma
\ ,\eeq
where  $k\in\N$ and $k$ are related to $\hb$ by   
 \beq
\label{hb0}
\hb
=\pm\frac{\pi}{2k+3}\ .\eeq
Define also
\[
[R+S]_{\tau\chi (T<0)}=F_{\hbar}(T,\tau\chi(T<0))^*\; S\;
F_{\hbar}(T,\tau\chi(T<0))\]
and
\[\tilde{\sigma}=F_{\hbar}(T,\tau\chi(T<0))^*\;\sigma\;
F_{\hbar}(T,\tau\chi(T<0))\ .\]
Since $F_{\hbar}(T,\tau\chi(T<0))$ is an unitary operator, 
 we see that   $([R+S]_{\tau\chi (T<0)},\tilde{\sigma})\in N_{\cal H}$. 
Therefore an operation 
\[\mlot_N:N^2_{\cal H}\longrightarrow N_{\cal H}\]
\[(R,\rho)\mlot_N (S,\sigma)=([R+S]_{\tau\chi(T<0)},
\tilde{\sigma})\ \]
is well defined. Moreover, it is not difficult to prove that this 
operation is associative. The operator domain $N$ together with 
the operation $\mlot$ 
is a braided quantum group.

Let $A$ be the operator domain as in  Example
\ref{pr:dziedzA}. 
Define an  operator map  
$$\varphi :N_{\cH}\rightarrow A_{\cH}$$
for any $(R,\rho)\in N_{\cH}$ by 
\beq
\label{a>n}
\varphi(R,\rho)=(R,\rho \chi(R<0))\ .\eeq 
The operator domain   $A^2$ and the group operation 
  $\mlot_A$ on $A$ were described in \cite{qexp}.
We prove now a  colorally, which enables us to apply  
 results obtained for $A$  in that paper to  \rf{rkN}.
We prove (see Corollary  \ref{n>a} below) that  
 for any $((R,\rho),(S,\sigma))\in N^2_{\cH}$
we have
\[\varphi((R,\rho)\mlot_N(S,\sigma))=\varphi(R,\rho)\mlot_A \varphi(S,\sigma)\ .\]
Moreover, next Colorally \ref{st:oexp} states that
\[
\fh(\varphi((R,\rho)\mlot_N (S,\sigma)))=\fh(\varphi(R,\rho))\fh(\varphi(S,\sigma))
\ .
\] 
\bfa 
\label{n>a}
Let $((R,\rho),(S,\sigma))\in N^2_{\cH}$.
Define
\[\hat{\tau} =(-1)^k\rho\chi(R<0)\sigma\chi(S<0)+
(-1)^k\sigma\chi(S<0)\rho\chi(R<0)\]
and
\[\tilde{\hat{\sigma}}=
\tilde{\sigma}\chi([R+S]
_{\tau\chi(T<0)}<0).\]
Then
\[
\tau\chi(T<0)=\hat{\tau}\mi
\tilde{\sigma}\chi(S<0)=\tilde{\hat{\sigma}}
\ .\]
\efa
\dow  Since $R$ ans $S$ satisfy Zakrzewski relation, it follows 
that $R$ commutes with  
$\sign S$ and $S$  commutes with 
$\sign R$. Hence
\[R\sigma\chi(S<0)=\sigma\chi(S<0)R\mi
S\rho \chi(R<0)=\rho \chi(R<0)S \ .\]
This means that if $((R,\rho),(S,\sigma))\in N_{\cal H}^2$ then 
  $(R,\rho\chi(R<0),S,\sigma\chi(S<0))$ 
satisfies assumptions of Theorem 6.1  \cite{qexp}.
By this theorem the  
sum $R+S$ has a selfadjoint extension, determined uniquely 
by a reflection  operator  
$\hat{\tau}$ such that 
\[\hat{\tau} =(-1)^k\rho\chi(R<0)\sigma\chi(S<0)+
(-1)^k\sigma\chi(S<0)\rho\chi(R<0)\ .\]
Since $\sigma$ anticommutes with  $R$, it follows  that  
$\chi(R<0)\sigma=\sigma\chi(R>0)$ and analogously for $\rho$ and $S$.
 Hence
\[\hat{\tau} =(-1)^k\rho\sigma\left\{\chi(R>0)\chi(S<0)+
\chi(S>0)\chi(R<0)\right\}=(-1)^k\rho\sigma\chi(\eh S^{-1}R<0)\ .\]
Comparing this result  with formula \rf{ttau} we see that 
$\hat{\tau}=\tau$.
It remains to prove that
\[\tilde{\hat{\sigma}}=
\tilde{\sigma}\chi([R+S]
_{\tau\chi(T<0)}<0).\]
Compute
\[\tilde{\hat{\sigma}}=\tilde{\hat{\sigma}}
\fh(T,\tau\chi(T<0))^*\sigma \chi(S<0)
\fh(T,\tau\chi(T<0))=\]
\[=\tilde{\sigma}
\fh(T,\tau\chi(T<0))^*\chi(S<0)
\fh(T,\tau\chi(T<0))=\tilde{\sigma}
\chi([R+S]_{\tau\chi(T<0)}<0)\;\;\od .\]\mqed
By  Theorem 6.1  \cite{qexp} we obtain 
\bst
\label{st:oexp}
Let $\left((R,\rho),(S,\sigma)\right)\in N_{\cal H}^2$ 
 and  let  $\ker S=\{0\}$ . 
%and 
%let  $\hat{T}=e^{i\frac{\hbar}{2}}S^{-1}R$ and $\hat{\tau}=(-1)^k\rho\hat{\sigma}+
%\overline{(-1)^k}\hat{\sigma}\rho$, gdzie 
%$(-1)^k =ie^{\frac{i\pi^2}{2\hb}}$.
Then 
\[\fh([R+S]_{\tau\chi (T<0)},\tilde{\sigma}\chi ([R+S]_{\tau\chi (T<0)}<0))=\]
\beq=
\fh(R,\rho\chi (R<0))\fh(S,\sigma\chi (S<0))
\label{oexp}
\ .\eeq
\est
The above Proposition explains why  $\fh$ is called the 
quantum exponential function. 
Moreover, the 
quantum exponential function  $\fh$, as  the classical exponential one,
 is the only one function (up to a parameter) 
satisfying the exponential equation  
\rf{oexp} (see theorem below).
\btw
\label{twstare}
Let $ \left((R,\rho),(S,\sigma)\right)\in N_{\cal H}^2 $ 
 and let $f\;:\;\Delta_{real}\rightarrow S^{1}$ be 
a measurable function.
Then the following  conditions are equivalent
\[\mbox{a).}\od\label{a:twstare}
%f(R,\rho\chi (R<0))f(S,\sigma\chi (S<0))=
%f([R+S]_{\tau\chi (T<0)},\tilde{\sigma}\chi ([R+S]_{\tau\chi (T<0)}<0))
f(\varphi(R,\rho))f(\varphi(S,\sigma))=
f(\varphi((R,\rho)\mlot_N (S,\sigma)))
\]
\[
%\[\!\!\!\!\!\!\!\!\!\!\!\!\!\!\!\!\!\\!\!\!
\mbox{b).}\od
 \mbox{ There exists  }
\mbox{$M\geq 0$ and $\mu=\pm 1$, }
\mbox{ such that  }\]
\[
\mbox{  }
f(\varphi(r,\rho))=\fh(\varphi(Mr,\mu  \rho))
\mbox{ for a.a.} (r,\rho)\in\R\times\{-1,1\}.\]    
\etw
\dow $\mbox{b). }\Longrightarrow \mbox{a). }$ 
We first consider the case  $M=0$.
Then
$$\fh(M r,\mu\rho)=\fh(0,\mu\rho)=1\ ,$$
because by   Theorem  1.1  \cite{qexp}
\[\lim _{r\rightarrow 0} \fh (r,\rho)=1\ .\]
It is easily seen that  if 
$M>0$ and $\mu=\pm 1$ then  
$$ \left((MR,\mu\rho),(MS,\mu\sigma)\right)\in N_{\cal H}^2 $$  
and $\ker MS=\{0\}$. Thus asumptions of Corollary
\ref{st:oexp} are satisfied and therefore   function 
$$f(r,\rho\chi(r<0))=\fh(Mr,\mu  \rho\chi(Mr<0))$$
 satisfies \rf{a:twstare}.

$\mbox{b). }\Longleftarrow \mbox{a). }$ 
By Corollary \ref{n>a} we may apply Theorem 7.1  \cite{qexp}, 
and b) follows.\hfill\qed
\subsection{The matrix representation of $N$}
Consider $((R,\rho),(S,\sigma))\in N_{\cal H}^2$.
Since $R\za S$, 
 operators  $R$ and $S$ commute with  
$\sign\; R$ and   $\sign\; S$.
Therefore we may introduce notation
\[\ba{lll}
 {\cH}_{++}&=&{\cH}(R>0)\cap {\cH}(S>0)\\
 {\cH}_{+-}&=&{\cH}(R>0)\cap {\cH}(S<0)\\ 
 {\cH}_{-+}&=&{\cH}(R<0)\cap {\cH}(S>0)\\
 {\cH}_{--}&=&{\cH}(R<0)\cap {\cH}(S<0)\\
\ea\]
Then  
${\cH}={\cH}_{++}\oplus {\cH}_{+-}\oplus {\cH}_{-+}\oplus {\cH}_{--}$.

Any vector  $\psi$ from the  space ${\cH}$ is represented by 
\[\psi=
\left[
\ba{ccc}
\psi_{++}\\
\psi_{+-}\\
\psi_{-+}\\
\psi_{--}
\ea\right]\]
where $\psi_{++}\in {\cH}_{++}$, $\psi_{+-}\in {\cH}_{+-}$, 
$\psi_{-+}\in {\cH}_{-+} $ and $\psi_{--}\in {\cH}_{--} $.
Therefore  operators acting on ${\cH}$ are represented by 
 $4\times 4$ matrices.
Moreover, since $\rho$ is selfadjoint and $\rho^2=
\chi(R\neq0)$, we see that maps
$\rho : {\cH}_{--}\rightarrow {\cH}_{-+} \mand  \rho : 
{\cH}_{-+}\rightarrow {\cH}_{--}$
are mutually inverse.
Similarly maps
$\rho : {\cH}_{+-}\rightarrow {\cH}_{++} \mand  \rho : 
{\cH}_{++}\rightarrow {\cH}_{+-}\ $
are mutually inverse.
Since $      \sigma$ is selfadjoint and $     \sigma ^2=
\chi(S\neq 0)$, we see that maps
$     \sigma   : {\cH}_{--}\rightarrow {\cH}_{+-} \mand       \sigma   : 
{\cH}_{+-}\rightarrow {\cH}_{--}$
are mutually inverse.
Also maps
$\si  : {\cH}_{-+}\rightarrow {\cH}_{++} \mand  \si  : 
{\cH}_{++}\rightarrow {\cH}_{-+}$
are mutually inverse.
Therefore  Hilbert spaces 
${\cH}_{-+},{\cH}_{+-}$ and ${\cH}_{--}$ and ${\cH}_{++}$ are 
unitarily equivalent. 
In what follows we simply assume that 
 ${\cH}_{++}={\cH}_{-+}={\cH}_{+-}={\cH}_{--}$ and 
 denote this  
Hilbert space by ${\cH}_+$. 

Then
\[
\rho=\left[\begin{array}{cccc}
0&I&0&0\\I&0&0&0\\0&0&0&I\\0&0&I&0
\end{array}\right],\hspace{20mm}\mand 
\sigma=\left[\begin{array}{cccc}
0&0&I&0\\0&0&0&I\\I&0&0&0\\0&I&0&0
\end{array}\right].
\]
Hence the  matrix representation of  operator $\tau:=(-1)^k\rho\sigma$ 
 is
\[
\tau
=
\left[\begin{array}{cccc}
0&0&0&I\\0&0&I&0\\0&I&0&0\\I&0&0&0
\end{array}\right].
\]

Since operators $R$ and $S$ commute with  $\sign\; R$ and   
with $\sign\; S$, they are represented 
by  
diagonal matrices
\[R=\left[\begin{array}{cccc}
R_+&0&0&0\\0&R_+&0&0\\0&0&-R_+&0\\0&0&0&\hspace{-2mm}-R_+
\end{array}\right],\hspace{10mm}\mand
S=\left[\begin{array}{cccc}
S_+&0&0&0\\0&-S_+&0&0\\0&0&S_+&0\\0&0&0&-S_+
\end{array}\right],
\]
where  $R_+$ and $S_+$ are restrictions to ${\cH}_{+}$
 of  $R$ and  $S$, respectively.
Clearly $R_+$ and  $S_+$  are selfadjoint and 
strictly positive  and  $R_+\za S_+$.
Moreover
\[
T=\eh S^{-1}R=\left[\begin{array}{cccc}
T_+&0&0&0\\0&-T_+&0&0\\0&0&-T_+&0\\0&0&0&T_+
\end{array}\right],
\]
where  $T_+=\eh S_+^{-1}R_+$.  $T_+$  
 is  selfadjoint  and 
strictly positive. 

%
%Operators $K$ and $S$ do commute with $\sign K$ and $\sign S$. Therefore they are represented by
%diagonal matrices. Kemembering that $\rho$ commutes with $K$ and $\hat{\sigma}$ commutes with $S$ we
%obtain 
%
%\subsection{Lemat }
\subsection{Matrix elements}
\label{subsect:uem}
%Let $\cK$ b/edzie spacei/a Hilberta.
%Let  $f\;:\;\Delta_{real}\rightarrow {\rm Unit} (\cK)$ b/edzie 
%funkcj/a measurableowsk/a and let for anyych $((R,\rho),(S,\sigma))
%\in N^2_{\cH}$ zachodzi
%\[f((R,\rho)\mlot_N(S,\sigma))=f(R,\rho)f(S,\sigma)\ .\]
%Then 
%\[
% f(r,\rho)f(s,\sigma)=f(s,\sigma)f(r,\rho)
%\]
%for anyych  $(r,\rho),(s,\sigma)\in
%\Delta_{real}$.
%\ele
%\dow
Consider strictly positive  operators $R$ and $S$, such that 
 $R\za S$. For example, $R_+$ and $S_+$ are such operators.
Since $R\za S$, by Corollary \ref{vN2} 
 the pair $(R,S)$ is 
unitarily equivalent to $(u\te e^{\po},v\te e^{\qo})$, 
where $u,v$ are unitary, 
selfadjoint and commuting operators, i.e. 
$\spe \;u,\spe\; v\subset\{-1,1\}$.
Assume that
\[R_o= e^{\po}\mi S_o= e^{\qo}\ .\]
%where 
%\[(\qo f)(x)=x f(x)\mand (\po f)(x)=\frac{\hbar}{i}\frac{d f(x)}{dx}\ .\]
%dzia/laj/acymi w  Hilbert space $\cK$.
%Wiemy, /ze 
%\[
% \hat{R}=\hat{R}^*\mand S=S^*\mand \hat{R}\za S \mand \hat{T}=e^{i\frac{\hbar}{2}}S^{-1}\hat{R}
%\label{pod}
%\]
Denote the complex conjugation operator  by  $J_o$. 
Then for any $w\in \Lkr$ 
\[ (J_o w)(t) =\ov{w(t)}\ ,\]
where  $t\in\R$.
Note that $J_o$
  is an antilinear operator 
   such that
\[
 (J_o)^2=I\mi (J_o)^*=J_o\mi J_o R_oJ_o=J_o  e^{\po}J_o=e^{-\po}=R_o^{-1}\]
and  
\[J_o S_o J_o=J_o  e^{\qo}J_o=e^{\qo}=S_o\ .
\]
Therefore  by Corollary  \ref{vN2} 
 for any $R\za S$ 
 there exists   an antilinear operator  $J$, 
  such that
\[ J^2=I\mi J^*=J \]
and
\beq
 JRJ=R^{-1}
\label{defJr}\eeq
and 
\beq JSJ=S
\label{defJs}
\ .\eeq
Since $J$ is  antilinear and 
$J^*=J$ then for any $w,v\in \cH$
\beq
\label{J*}
\langle w|Jv\rangle=\langle v|J^*w\rangle=\langle v|Jw\rangle\ .
\eeq
Moreover, define an  operator $F$ by
\beq
F=e^{i\frac{\pi}{4}}e^{-i\frac{\log^2 S}{\hbar}} 
e^{-i\frac{\log^2 T}{2\hbar}}
\label{defF}
\ .\eeq
Note that
\beq F^*=F^{-1}
\label{fg}\ .\eeq
It is not very difficult to see that if $R=R_o$ and $S=S_o$ then $F$ is the Fourier transform. 

By  \rf{defJr}, \rf{defJs} and  (\ref{defF}) 
\beq
\label{fr}
 FRF^{-1}=S\eeq	
and
\[  
%\label{fsuem}
FSF^{-1}=R^{-1}\ .\]
Moreover
\beq
\label{fj} 
 FJ=JF^{-1}\mand F^{-1}J=JF
\eeq
and 
\beq
\label{ftf}
 JTJ=FT^{-1}F^{-1}
.\eeq	
%\bpr
%if przyjmiemy \[(\qo f)(x)=x f(x)\mand (\po f)(x)=i\hbar\frac{d f(x)}{dx}\]
%\[\hat{R}=e^{-\po}\mand S=e^{\qo}\]
%where 
%\[(\qo f)(x)=x f(x)\mand (\po f)(x)=i\hbar\frac{d f(x)}{dx}\]
%to then 
% \[ \hat{T}=e^{i\frac{\hbar}{2}}S^{-1}R=e^{-\po-\qo}\]
%operator $J$ to sprz/e/zenie zespolone, za/s $F$ is transformat/a 
%Fouriera.
%Zwr/o/cmy jednak uwag/e, /ze powy/zsze wzory wyprowadzone are tylko w 
%oparciu o zale/zno/sci (\ref{pod}),  (\ref{defJ}) and  (\ref{defF}), 
%i are prawdziwe niezale/znie od reprezentacji.
%\label{np2}
%\epr
We use the notion of generalized eigenvectors.
It is well known that a selfadjoint operator with continuous spectrum 
 acting on $\cH$ does not have eigenvectors. Still one can show that in the general case 
the generalized eigenvectors are continiuous linear functionals 
on a certain dense locally convex subspace $\Phi\subset \cH$, provided with a much 
stronger topology than $\cH$. Then we get the same formulae as
 for discreet spectrum provided we replace scalar product by the duality relation  between $\Phi$ and $\Phi^\prime$.   
 This will be explainded 
 by the example below, for general considerations see \cite{maueig}.

%Antyliniowy funkcjona/l  $x_a$ na $W$ nazywamy {\em uog/olnionym wektorem w/lasnym} operatora $A$, odpowiadaj/acym warto/sci w/lasnej $a$, gdy dla 
%dowolnych $w\in W$ zachodzi
%\[
%\label{defuww}
%\langle Aw|x\rangle=a \langle w|x\rangle\ .\]
%W dalszej cz/e/sci b/edziemy pisa/c dla uproszczenia
% \[A|x_a\rangle=\lambda |x_a\rangle\]
%zamiast \rf{defuww}.

\bpr 
Let $\cH=\Lkr$
% $W=\sz$, 
and
\beq
\label{np2}
 R=e^{\po} \mi S=e^{\qo} \mi T=\eh S^{-1}R=e^{\po-\qo}\ .\eeq
These operators have continuous spectra, so they 
do not have eigenvectors.
%, i.e. there are no functions $\Omega_r$, 
%$\Phi_s$ and $\Psi_t\in\Lkr$ such that 
 There are however tempered distributions on $\R$ 
such that for every function $f$
  from the  Schwartz space of smooth functions on $\R$ 
 decreasing rapidly at infinity  $\sz$ we have 
\beq
\langle f|R|\Omega_r\rangle=r\langle f|\Omega_r\rangle\mi 
\langle f|S|\Phi_s\rangle=s\langle f|\Phi_s \rangle  \mi 
\langle f |T|\Psi_t\rangle =t\langle f|\Psi_t\rangle . 
\label{eigenvectors}\eeq
Such  $|\Omega_r\rangle$, $|\Phi_s\rangle$ and $|\Psi_t\rangle$ 
 are called generalized eigenvectors of 
 operators  $R$,$S$ and $T$ with  eigenvalues respectively $r$, $s$ and $t$.

An example of  generalized eigenvectors of
 operators \rf{np2} is
\[|\Omega_r\rangle=\frac{1}{\sqrt{2\pi\hbar}}e^{\frac{i}{\hbar}x\log r}\mi
|\Phi_s\rangle=\delta(\log s-x)\mi
|\Psi_t\rangle=\frac{1}{\sqrt{2\pi\hbar}}e^{i\frac{x^2}{2\hbar}}
e^{i\frac{x\log t}{\hbar}}\ .\]

Moreover, we will use notation of a type 
$\langle \Omega_r|\Phi_s\rangle$.

It should be understood in the following way: for any $f\in \sz$ we have
\beq\langle \Omega_r|f\rangle=\int_{\R} \langle \Omega_r|\Phi_s\rangle 
\langle \Phi_s|f\rangle ds\ .
\label{dlugie}\eeq
To shorten notation from now on we skip  
 the integration symbol, i.e. we write 
\[\langle \Omega_r|f\rangle=\langle \Omega_r|\Phi_s\rangle 
\langle \Phi_s|f\rangle  .\]
instead of \rf{dlugie}
The generalized eigenvectors $\Omega$ are said to have     
 the  {\em Dirac $\delta$ normalization} if  
\[\langle \Omega_r|\Omega_s\rangle =\delta (r-s)\ ,\]
where $\delta$ is the Dirac  $\delta$ distribution.
Note that generalized eigenvectors $\Omega$, $\Phi$ and $\Psi$ given above 
have  the   Dirac $\delta$ normalization.
\epr

Let $|\Omega_r\rangle$ be a generalized eigenvector of $R$  
with real eigenvalue $r$ and with  Dirac delta normalization.
Analogously, let  $|\Phi_s\rangle$ and $|\Psi_t\rangle$ denote 
generalized eigenvectors of 
 operators  $S$ and $T$ with real eigenvalues respectively  $s$ and $t$
 and with Dirac delta normalization.
%\[R|\Omega_r\rangle=r|\Omega_r\rangle\mi S|\Phi_s\rangle=s|\Phi_s\rangle \mi T|\Psi_t\rangle=t|\Psi_t\rangle\ .\]
%%Zwr/o/cmy uwag/e, /ze $r,s,t\in \R$.
%
%hence na mocy cytowanego na samym pocz/atku dowodu  swierdzenia \ref{vN2}
% tak/ze dla og/olnych $R\za S$ mo/zna 
%  dobra/c  uog/olnione, unormowane do delty Diraca,  wektory w/lasne 
% operator/ow $R$, $S$ and $T$, tak, aby zachodzi/lo \ref{dobranie}.
Note that by (\ref{defJs})
% and (\ref{fj}) 
 for any $w\in \cH$ 
\[ \langle w|S\Phi_s\rangle=\langle w|JSJ\Phi_s\rangle\ .\]
Since $J$ is an antilinear operator and $S$ commutes with $J$, it 
follows
\[\langle w|JSJ\Phi_s\rangle=\langle SJ\Phi_s|J w\rangle=
\langle JS\Phi_s|J w\rangle=\langle w|J^2S\Phi_s\rangle=
%Z drugiej strony, poniewa/z $\Phi_s$ is uog/olnionym wektorem w/lasnym 
% $S$ to
\langle w|S\Phi_s\rangle=s\langle w|\Phi_s\rangle\ .\]
Comparing above formulae we obtain
\[\langle SJ\Phi_s| J w\rangle=s\langle w|\Phi_s\rangle
=s\langle w|J^2\Phi_s\rangle
=s\langle J\Phi_s|Jw\rangle
\ ,\]
so
%, poniewa/z $s\in \R$,
\[S|J\Phi_s\rangle=s|J\Phi_s\rangle\ .\]
Hence
\beq\label{uem1}
|J\Phi_s\rangle=|\Phi_s\rangle\ .\eeq
Moreover observe, that by \rf{fr} 
\[ FRF^{-1}|\Phi_r\rangle =r |\Phi_r\rangle\ .\]
Hence
\[ R|F^{-1}\Phi_r\rangle =r |F^{-1}\Phi_r\rangle\ ,\]
so
\beq \label{uem2} |F^{-1}\Phi_r\rangle = |\Omega_r\rangle\ .\eeq
Applying \rf{uem1} and \rf{uem2}, and then\rf{fg} and \rf{fj},
 and again  \rf{uem2} and \rf{J*} and \rf{uem1}, we obtain 
\[\langle\Omega_r|\Phi_s \rangle=\langle F^{-1}\Phi_r|J \Phi_s \rangle=\langle\Phi_r|F\;J \Phi_s \rangle=\]
\[=\langle\Phi_r|J\;F^{-1}\; \Phi_s \rangle=\langle\Phi_r|J\;\Omega_s\rangle=\langle\Omega_s|J\;\Phi_r\rangle
=\langle\Omega_s|\;\Phi_r\rangle ,\]
so
\beq
\label{omfi}
\langle\Omega_r|\Phi_s \rangle=\langle\Omega_s|\;\Phi_r\rangle
\ .\eeq
We proceed to deriving our next formula.
Note that
\[
e^{i\frac{\log^2 t}{2\hb}} \langle \Omega_r|\Psi_t\rangle
\langle\Psi_t|            \Phi_s \rangle=
\langle \Omega_r|e^{i\frac{\log^2 T}{2\hb}}\Psi_t\rangle
\langle\Psi_t|            \Phi_s \rangle
=
\langle F^{-1}\Phi_r|e^{i\frac{\log^2 T}{2\hb}}\Psi_t\rangle
\langle\Psi_t|            \Phi_s \rangle=\]
\[=
\langle \Phi_r|e^{i\frac{\pi}{4}}e^{-i\frac{\log^2 S}{\hbar}} 
e^{-i\frac{\log^2 T}{2\hbar}}e^{i\frac{\log^2 T}{2\hb}}\Psi_t\rangle
\langle\Psi_t|            \Phi_s \rangle=
\langle \Phi_r|e^{i\frac{\pi}{4}}e^{-i\frac{\log^2 S}{\hbar}} \Psi_t\rangle
\langle\Psi_t|            \Phi_s \rangle=\]
\[=
e^{i\frac{\pi}{4}}
\langle e^{i\frac{\log^2 S}{\hbar}}\Phi_r| \Psi_t\rangle
\langle\Psi_t|            \Phi_s \rangle
=e^{i\frac{\pi}{4}}e^{-i\frac{\log^2 r}{\hbar}}  
\langle \Phi_r| \Psi_t\rangle \langle \Psi_t|   \Phi_s \rangle
\ .\]
It means, we have the formula
\beq
\label{pom1}
e^{-i\frac{\pi}{4}}e^{i\frac{\log^2 t}{2\hb}} \langle \Omega_r|\Psi_t\rangle
\langle\Psi_t|            \Phi_s \rangle
=e^{-i\frac{\log^2 r}{\hbar}}  
\langle \Phi_r| \Psi_t\rangle \langle \Psi_t|  \Phi_s \rangle
\ .\eeq
We prove now that
\beq
\label{ef}
\langle\Omega_r|\vt(\log T)^*|\Phi_{s}\rangle=c_\hbar^\prime  
e^{-i\frac{\log^2 s}{\hbar}}\langle\Phi_{s}|\vt(\log T)|\Phi_r\rangle
\ .\eeq
Compute the left hand side of this formula
 \[
LHS=
\langle\Omega_r|\vt(\log T)^*|\Phi_{s}\rangle=
\langle\Omega_r|\vt(\log T)^*|\Psi_{t}\rangle
\langle\Psi_{t}|\Phi_{s}\rangle=\]
\[\label{L}=\ov{\vt(\log t)}
\langle\Omega_r|\Psi_{t}\rangle
\langle\Psi_{t}|\Phi_{s}\rangle\ .\]
Moreover by formula 1.36  \cite{qexp}, for any  $t\in\R$
\beq
\label{vts}
\overline{\vt(\log t)}=e^{-i\frac{\pi}{4}}c_\hbar^\prime e^{i\frac{\log^2 t}{2\hb}}
\vt(-\log t)
\ ,\eeq
where $ c_\hbar^\prime=
e^{i(\frac{\pi}{4}+\frac{\hb}{24}+\frac{\pi^2}{6\hb})}$. 
Hence
\[L=e^{-i\frac{\pi}{4}}c_\hbar^\prime e^{i\frac{\log^2 t}{2\hb}}
\vt(-\log t)
\langle\Omega_r|\Psi_{t}\rangle
\langle\Psi_{t}|\Phi_{s}\rangle\ .\]
Compute now the right hand side of \rf{ef}
\beq
\label{P}
RHS=
c_\hbar^\prime  
e^{-i\frac{\log^2 s}{\hbar}}\langle\Phi_{s}|\vt(\log t)|\Psi_{t}\rangle
\langle\Psi_{t}|\Phi_r\rangle
=c_\hbar^\prime  \vt(\log t)
e^{-i\frac{\log^2 s}{\hbar}}\langle\Phi_{s}|\Psi_{t}\rangle
\langle\Psi_{t}|\Phi_r\rangle
\ .\eeq	
Note that by \rf{defF} and \rf{fg}
\[T|e^{i\frac{\log^2 S}{\hbar}}J\Psi_{t}\rangle=
e^{-i\frac{\log^2 T}{2\hbar}}
e^{i\frac{\pi}{4}}TF^{-1}|J\Psi_{t}\rangle.\]
Moreover by \rf{ftf}
\[e^{-i\frac{\log^2 T}{2\hbar}}
e^{i\frac{\pi}{4}}TF^{-1}J|\Psi_{t}\rangle
=e^{-i\frac{\log^2 T}{2\hbar}}
e^{i\frac{\pi}{4}}F^{-1}JT^{-1}JJ|\Psi_{t}\rangle=\]
\[=t^{-1}e^{-i\frac{\log^2 T}{2\hbar}}
e^{i\frac{\pi}{4}}F^{-1}J|\Psi_{t}\rangle=
t^{-1}|e^{i\frac{\log^2 S}{\hbar}}J\Psi_{t}\rangle\ ,\]
so
\[|e^{i\frac{\log^2 S}{\hbar}}J\Psi_{t}\rangle
=|\Psi_{t^{-1}}\rangle\ .\]
Hence
\[ \vt(\log t)
\langle\Phi_{s}|\Psi_{t}\rangle
\langle\Psi_{t}|\Phi_r\rangle
=\vt(-\log t)
\langle\Phi_{s}|\Psi_{t^{-1}}\rangle
\langle\Psi_{t^{-1}}|\Phi_r\rangle
=\]
\[=\vt(-\log t)
\langle\Phi_{s}|e^{i\frac{\log^2 S}{\hbar}}J\Psi_{t}\rangle
\langle e^{i\frac{\log^2 S}{\hbar}}J\Psi_{t}|\Phi_r\rangle
=\]
\[=\vt(-\log t)e^{-i\frac{\log^2 r}{\hbar}}
e^{ i\frac{\log^2 s}{\hbar}}
\langle\Phi_{s}|J\Psi_{t}\rangle
\langle J\Psi_{t}|\Phi_r\rangle
=\]
\[=\vt(-\log t)e^{-i\frac{\log^2 r}{\hbar}}
e^{ i\frac{\log^2 s}{\hbar}}\langle  \Psi_{t}|J\Phi_s\rangle
\langle J\Phi_  r|\Psi_{t}\rangle
=\]
\[=\vt(-\log t)e^{-i\frac{\log^2 r}{\hbar}}
e^{ i\frac{\log^2 s}{\hbar}}
\langle \Psi_{t}|\Phi_s\rangle
\langle\Phi_{r}|\Psi_{t}\rangle
\ .
\]
Therefore
\beq
\label{pom20}
\vt(\log t)
\langle\Phi_{s}|\Psi_{t}\rangle
=\vt(-\log t)e^{-i\frac{\log^2 r}{\hbar}}
e^{ i\frac{\log^2 s}{\hbar}}
\langle \Psi_{t}|\Phi_s\rangle
\langle\Phi_{r}|\Psi_{t}\rangle
\ .
\eeq
In fact, we have even proved a more general formula, namely
for any measurable function  $f$ we have
\beq
\label{pom2}
f(\log t)
\langle\Phi_{s}|\Psi_{t}\rangle
\langle\Psi_{t}|\Phi_r\rangle
=
f(-\log t)e^{-i\frac{\log^2 r}{\hbar}}
e^{ i\frac{\log^2 s}{\hbar}}
\langle \Psi_{t}|\Phi_s\rangle
\langle\Phi_{r}|\Psi_{t}\rangle
\ .
\eeq
We use this formula in our forthcomming paper \cite{paper4} on the quantum
 'az+b' group at roots of unity.

Comparing \rf{P} and \rf{L} and using \rf{pom1} and \rf{pom20} 
we get \rf{ef}.
Using again formulae (1.36) \cite{qexp} with   
  $z=-\log t-i\pi$, where $t\in\R$, and 
  \rf{hb0}, we   obtain
\[\overline{\vt(\log t-i\pi)}=c_\hbar^\prime (-1)^k e^{-i\frac{\pi}{4}}
e^{i\frac{\overline{\log t} ^2 }{2\hb}} e^{-\frac{\pi}{\hb}\overline{\log t}}
\vt(-\overline{\log t}-i\pi) \ .\]
Using the above formula  and the same method as  
for derivation of \rf{ef}, we obtain formulae for some matrix elements 
we will soon find very useful
\beq
\label{efw}
\langle\Omega_r|\vt(\log T-i\pi I)^*|\Phi_{\tilde{s}}\rangle=i
(-1)^k c_\hbar^\prime  
e^{-i\frac{\log^2 \tilde{s}}{\hbar}}\langle\Phi_{\tilde{s}}|
T^\frac{\pi}{\hb}\vt(\log T-i\pi I)|\Phi_r\rangle
\eeq
and
\beq
\label{efwt}
\langle\Omega_r|T^\frac{\pi}{\hb}\vt(\log T-i\pi I)^*
|\Phi_{\tilde{s}}\rangle=i(-1)^k c_\hbar^\prime  
e^{-i\frac{\log^2 \tilde{s}}{\hbar}}\langle\Phi_{\tilde{s}}|
\vt(\log T-i\pi I)|\Phi_r\rangle
\ .\eeq 
%
%********************
% 
\subsection{Unitary representations of $N$}
We find now a formula for all   
 unitary reprezentations of    $N$ acting on a Hilbert space $\cK$. 
\bde
\label{repuni}
 Let for any Hilbert space $\cH$ there exists a map
  $$V_{\cH}:N_{\cal H}\longrightarrow 
\mbox{\rm Unit}(\cK\te \cH)$$ 
such that
\bit
\item[1.]{For any  $(R,\rho)\in N_{\cH}$
 and for any operators $ v\in \mbox{\rm Unit}(\cH,\cH ^\prime)$ 
we have
$$(\id _{\cK}\te v^*)V_{\cH}(R,\rho) (\id _{\cK}\te v)= V_{\cH^\prime}
(v^*Rv,v^*\rho v)\ .$$ }
\item[2.]{For any space with measure $(\Lambda,\mu)$ and for any 
 measurable field of  Hilbert spaces $\{\cH(\lambda)\}_{\lambda\in\Lambda}$ 
and for any  measurable fields of closed operators  
$\{R(\lambda)\}_{\lambda\in\Lambda}$ and $\{\rho(\lambda)\}_{\lambda\in\Lambda}$,
 such that
$(R(\lambda),\rho(\lambda))\in N_{\cH(\lambda)}$
we have
$$\int _{\Lambda}^\oplus V_{\cH(\lambda)}(R(\lambda) ,\rho (\lambda)) d\mu (\lambda)
=V_{\int _{\Lambda}^\oplus \cH(\lambda)d\mu (\lambda)}( \int _{\Lambda}^\oplus R(\lambda)d\mu (\lambda),
\int _{\Lambda}^\oplus \rho(\lambda)d\mu (\lambda)) .$$}
%where $( \hat{R}_1\oplus \hat{R}_2,
% \hat{\rho}_1 \oplus \hat{\rho}_2)\in A_{{\cH}_1\oplus {\cH}_2}$. }
\item[3.]{For any $ 
\left((R,\rho),(S,\sigma)\right)\in N_{\cal H}^2$  we have 
\beq
 \label{repn} V_{\cH}(R,\rho) V_{\cH} (S,\sigma) =V_{\cH}\left(R,\rho)\mlot_{N} 
(S,\sigma)
\right)\ .\eeq
%U([\hat{R}+\hat{S}]_\phi,\tilde{\hat{S}})            $$
%where $\phi =\pm \hat{\rho}\\hat{S}igma \chi(\eh \hat{S} \hat{R}<0)$.  
  }
\eit
Then we call  $V$ a  {\bf  unitary representation} 
of $N$ on Hilbert space $\cK$.
\ede
In what follows we omit the subscript  $\cH$ 
in $V_{\cH}$. 

We prove now a formula for all 
  unitary representations of  $N$ on a Hilbert space $\cK$.
\btw
\label{twN}
A map  $$V_{\cK}:N_{\cal H}\longrightarrow 
\mbox{\rm Unit}(\cK\te \cH)$$ 
is a unitary representation of  $N$
iff, there exists $(M,\mu)\in N_{\cK}$, such that  
 for any
 \mbox{$(R,\rho)\in N_{\cH}$} we have
\[V(R,\rho)=\fh(M\te R,\mu \te \rho\chi(M\te R<0))\ .\]
\etw
By  $\fh(M\te R,\mu \te \rho\chi(M\te R<0))$ we mean 
\[\fh(M\te R,\mu \te \rho\chi(M\te R<0))=\int_{\R\times\{-1,1\}} 
\fh(R,\rho\chi(R<0))\te dE_{M,\mu}(z),\]
where $dE_{M,\mu}$ is the joint spectral measure of strongly 
 commuting  
operators  $M$ and $\mu$,  acting on Hilbert space $\cK$.

We proceed to prove  Theorem \ref{twN}.

\dow $\Leftarrow$ Observe that  if $(R,\rho)\in N_{\cal H}$
and   $(M,\mu)\in N_{\cal K}$ then also 
 $$(M\te R,
\mu\te\rho)\in N_{{\cal K}\te{\cal H}}\ .$$
%Ponadto, na mocy Stwierdzenia \ref{an2} then  
%$$(M\te R,
%(\mu\te\rho)\chi(M\te R<0))\in 
%A_{{\cal K}\te{\cal H}},$$
%za/s ze Stwierdzenia \ref{an3} wynika, /ze:
%$$(M\te R,\mu\te\rho)\mlot_N (M\te S,
%\mu\te\sigma)=$$
%$$=(M\te R,(\mu\te\rho)\chi(M\te R<0))
%\mlot_A (M\te S,
%(\mu\te\sigma)\chi(M\te S<0)).$$
Therefore we may apply  Theorem \ref{twstare}, which is our claim.

$\Rightarrow$ We follow the proof of 
  Theorem 4.2  \cite{oe2}
%, tzn. wyka/zemy, /ze 
%$ f(r,\rho)f(s,\sigma)=f(s,\hat{\sigma})f(r,\rho)$, dla $(r,\rho),(s,\hat{\sigma})\in
%\Delta_{real}$.
. We first outline the proof.
We show that if  $V$ is a unitary representation of $N$,
then 
%for any  $(R,\rho,S,\sigma)\in N_{\cal H}^2$  
\[V(r,\varrho)V(s,\sigma)=V(s,\sigma)V(r,\varrho)\]
for any $r,s\in \R\setminus\{0\}$ and $\varrho,\sigma\in \{-1,1\}$.
Then we find formula for all unitary representions 
 of  $N$ acting on  $\C$. Using spectral decomposition theorem completes the proof.

Our proof starts with the  observation that  
 since  $R$ and $\rho$ commute and 
$\spe \rho =\{-1,1\}$, the function $V$  may be written  as
\beq
\label{v}
V(R,\rho)=V_1(R)+(I_{\cal K}\te \rho) V_2(R),
\eeq
where
\[V_1(R)=\frac{1}{2}(V(R,1)+V(R,-1))\mi 
V_2(R)=\frac{1}{2}(V(R,1)-V(R,-1)).\]
Then
\[V(R,\rho)=\left[\ba{cccc}V_1(R_+)&V_2(R_+)&0&0\\
V_2(R_+)&V_1(R_+)&0&0\\
0&0&V_1(-R_+)&V_2(-R_+)\\0&0&V_2(-R_+)&V_1(-R_+)\ea\right]
\]
and
\[V(S,\sigma)=\left[\ba{cccc}V_1(S_+)&0&V_2( S_+) &0\\
0&V_1(-S_+)&0&V_2(-S_+)\\
V_2( S_+) &0&V_1(S_+)&0\\0&V_2(-S_+)&0&V_1(-S_+)\ea\right]
.\]
Hence
\[V(R,\rho)V(S,\sigma)=\]
\[=
\left[\ba{cccc}V_1( R_+)V_1( S_+)&V_2( R_+)V_1(-S_+)
&V_1( R_+)V_2( S_+)
&V_2( R_+)V_2(-S_+)\\
V_2( R_+)V_1( S_+)&V_1( R_+)V_1(-S_+)&V_2( R_+)
V_2( S_+)&V_1( R_+)V_2(-S_+)\\
V_1(-R_+)V_2( S_+)&V_2(-R_+)V_2(- S_+)&V_1(-R_+)
V_1( S_+)&V_2(-R_+)V_1(-S_+)\\
V_2(-R_+)V_2( S_+)&V_1(-R_+)V_2(-S_+)&V_2(-R_+)
V_1( S_+)&V_1(-R_+)V_1(-S_+)
\ea\right]
\ .\]
Moreover
\[\fh (T,\tau\chi(T<0))=\]
\[= 
\left[
\ba{cccc}
\vt (\log T_+)&0&0&0\\
0&\vt (\log T_+- i\pi I)&i(-1)^k T_+^{\frac{\pi}{\hb}}
\vt(\log T_+- i\pi I)&0\\
0&i(-1)^k T_+^{\frac{\pi}{\hb}}
\vt(\log T_+- i\pi I)&\vt (\log T_+- i\pi I)&0\\
0&0&0&\vt (\log T_+)\\
\ea\right]
\ .\]
Hence
\[X:=V([R+S]_{\tau\chi(T<0)},\tilde{\sigma})=
\left[
\ba{cccc}
X(1,1)&X(1,2)&X(1,3)&0\\
X(2,1)&X(2,2)&X(2,3)&X(2,4)\\
X(3,1)&X(3,2)&X(3,3)&X(3,4)\\
0&X(4,2)&X(4,3)&X(4,4)
\ea\right]\ ,\]
where
\[\ba{ccc}
X(1,1)&=&\vt(\log T_+)^*\;V_1(S_+  )\;\vt(\log T_+)\\
X(1,2)&=&i(-1)^k\vt (\log T_+)^*V_2(S_+)(T_+)^{\frac{\pi}{\hb}}
\vt(\log T_+- i\pi I)\\
X(1,3)&=&\vt (\log T_+)^*V_2(S_+)
\vt(\log T_+- i\pi I)\\
X(2,1)&=&-i(-1)^k(T_+)^{\frac{\pi}{\hb}}
\vt(\log T_+- i\pi I)^*V_2(S_+)\vt (\log T_+)\\
X(2,2)&=&
\vt(\log T_+-i\pi I)^*\;V_1(-S_+  )\;\vt(\log T_+-i\pi and )\od +\\
&+&  T_+^{\frac{\pi}{\hb}}\;
\vt(\log T_+-i\pi I)^*\;
 V_1(S_+)   \; T_+^{\frac{\pi}{\hb}}\;\vt(\log T_+-i\pi I) \\
X(2,3)&=&
i\;(-1)^k\;\vt(\log T_+-i\pi I)^*\;  V_1 (-S_+  )\;
 T_+^{\frac{\pi}{\hb}}\;\vt(\log T_+-i\pi I)
\od  -\\
& -&  i\;(-1)^k\; T_+^{\frac{\pi}{\hb}}\;
\vt(\log T_+-i\pi I)^*\;   V_1 (S_+)   \;\vt(\log T_+-i\pi I)\\
X(2,4)&=&\vt(\log T_+- i\pi I)^*V_2(-S_+) \vt (\log T_+)\\
X(3,1)&=&\vt(\log T_+- i\pi I)^*V_2(S_+)\vt (\log T_+)\\
X(3,2)&=& 
\; -i\;(-1)^k\; 
T_+^{\frac{\pi}{\hb}}\;\vt(\log T_+-i\pi I)^*\; V_1(-S_+  )\;
\vt(\log T_+-i\pi I)\od  +\\
 &+& 
i\;(-1)^k\;\vt(\log T_+-i\pi I)^*\;  V_1(S_+)   \;
  T_+^{\frac{\pi}{\hb}}\;\vt(\log T_+-i\pi I)\\
X(3,3)&=&  T_+^{\frac{\pi}{\hb}}\;
\vt(\log T_+-i\pi I)^*\;V_1(-S_+   )\;
 T_+^{\frac{\pi}{\hb}}\;\vt(\log T_+-i\pi I)\od  +\\
 &+&  
\vt(\log T_+-i\pi I)^*\;   V_1( S_+)   \;
\vt(\log T_+-i\pi I)\\
X(3,4)&=&-i(-1)^k(T_+)^{\frac{\pi}{\hb}}
\vt(\log T_+- i\pi I)^*V_2(-S_+)\vt (\log T_+)\\
X(4,2)&=&\vt (\log T_+)^*V_2(-S_+)\vt(\log T_+- i\pi I)\\
X(4,3)&=&i(-1)^k\vt (\log T_+)^*V_2(-S)(T_+)^{\frac{\pi}{\hb}}
\vt(\log T_+- i\pi I)\\
X(4,4)&=&\vt(\log T_+)^*\;V_1(-S_+  )\;\vt(\log T_+)
\ea\ .\]
We prove that
\bst
\label{przemN}For any $r,s\in\R$ and $\varrho,\sigma\in\{-1,1\}$ we have
\beq
\label{wz:przemN}
V(r,\varrho)V(s,\sigma)=V(s,\sigma)V(r,\varrho)
.\eeq
Moreover 
\[V_2(r)V_2(-s)=0=V_2(-s)V_2(r).\]
\est
\dow We prove that
\beq
\label{1++}
%1
V_1(r  )V_1(s  )=V_1(s  )V_1(r  )
\eeq
\beq
%2
\label{1--}
V_1(-r  )V_1(-s  )=V_1(-s  )V_1(-r  )
\eeq
\beq
%3
\label{1+1-}
V_1(r)V_1(-{s})=V_1(-{s})V_1(r)
\eeq
\beq
\label{2+2-}
V_2(r)V_2({s})= V_2({s})V_2({r})
\eeq
\beq
\label{2-2+}
V_2(-r)V_2(-{s})= V_2(-{s})V_2(-{r})
\eeq
\beq
\label{2--}
V_2(r)V_2(-s)=0=V_2(-s)V_2(r)
\eeq
\beq
\label{1-2+}
V_1(r)V_2({s})= V_2({s})V_1({r})
\eeq
\beq
\label{2-1+}
V_1(-r)V_2(-{s}) = V_2(-{s})V_1(-{r}) 
\eeq
\beq
\label{1+2-}
V_1({r})V_2(-{s})= V_2(-{s})V_1({r})
\eeq	
\beq
\label{1+2+}
V_1(-r)V_2({s})= V_2({s})V_1(-r)
\ .\eeq
The formulae  \rf{omfi}, \rf{ef} and 
\rf{efw} and \rf{efwt}  will be of great use throughout the proof.
First we prove the formula  \rf{1++}.

%\subsubsection{$V_1({r})V_1({s})=V_1({s})V_1({r})$}
Compute
\[\langle\Omega_{r}|V_1({R_+})V_1(S_+)|\Phi_{s}\rangle=
\langle\Omega_{r} |X(1,1)|\Phi_{s}\rangle=
\langle\Omega_{r}|\vt(\log T_+)^*\;V_1({S_+})\;\vt(\log T_+)|\Phi_{s}\rangle
 \ .\]
Hence
\[V_1({r})V_1({s})=\langle\Omega_{r}|\Phi_{s}\rangle^{-1}\langle\Omega_{r}|\vt(\log t)^*|\Psi_t\rangle
\langle\Psi_t|\Phi_{\tilde{{s}}}\rangle
V_1(\tilde{{s}})|\langle\Phi_{ \tilde{{s}}}
|\Psi_{\tilde{t}}\rangle\langle\Psi_{\tilde{t}}|\vt(\log t)|\Phi_s\rangle .\]
Therefore applying  \rf{ef} we get 
$$V_1({r})V_1({s})=V_1(\tilde{{s}})\langle\Omega_{r}|
\Phi_s\rangle^{-1} 
\langle\Omega_{r}|\vt(\log T_+)^*|\Phi_{\tilde{{s}}}\rangle 
\langle\Phi_{\tilde{{s}}}|\vt(\log T_+)|\Phi_{s}\rangle=$$ 
$$ =c_\hbar ^\prime e^{-i\frac{\log^2 \tilde{{s}}}{\hbar}} 
 V_1(\tilde{{s}})\langle\Omega_{r}|\Phi_{s}\rangle^{-1} 
\langle\Phi_{\tilde{{s}}}|\vt(\log T_+)|\Phi_{r}\rangle 
 \langle\Phi_{\tilde{{s}}}|\vt(\log T_+)|\Phi_{s}\rangle .$$ 
Note that  $\langle\Phi_{\tilde{{s}}}|\vt(\log T_+)|\Phi_{r}\rangle 
 \langle\Phi_{\tilde{{s}}}|\vt(\log T_+)|\Phi_{s}\rangle$ is 
symmetric with respect to swapping  
  ${r}$ and ${s}$. By (\ref{omfi}) the same holds for  
 \mbox{$\langle\Omega_{r}|\Phi_{s}\rangle$} and the remaining terms  
of the above formula depend on neither  ${r}$ nor  $s$.
 Therefore
$$V_1({r})V_1({s})=V_1({s})V_1({r}),$$
i.e.  \rf{1++} holds. 

In the same manner one can prove formulae  
\rf{1--} and \rf{2+2-} and \rf{2-2+}.
 
The proof of the formula  \rf{1+1-} is slightly different. 
Compute
%\subsubsection{$V_1({r})V_1(-{s})=V_1(-{s})V_1({r})$}
\[\langle\Omega_{r}|V_1(R_+)V_1(-S_+)|\Phi_s\rangle=\langle\Omega_{r}|
X(2,2)|\Phi_{s}\rangle=\]
\[=
\langle\Omega_{r}|\;\vt(\log T_+-i\pi I)^*\;V_1(-{S_+})\;
\vt(\log T_+-i\pi I)\;|\Phi_{s}\rangle+\]
 \[+\;\langle\Omega_{r}|
 T_+^{\frac{\pi}{\hb}}\;\vt(\log T_+-i\pi I)^*\;V_1({S_+})\;
  T_+^{\frac{\pi}{\hb}}\;\vt(\log T_+-i\pi I) |\Phi_{s}\rangle. \]
Hence
\[V_1({r})V_1(-{s})=\langle\Omega_{r}|\Phi_{s}\rangle^{-1}\;\times\]
\[\times\;\left\{\;\langle\Omega_{r}|\;
\vt(\log t-i\pi I)^*|\Psi_t\rangle
\langle\Psi_t|\Phi_{\tilde{{s}}}\rangle
V_1(-\tilde{{s}})
|\langle\Phi_{ \tilde{{s}}}|\Psi_{\tilde{t}}\rangle\langle\Psi_{\tilde{t}}|\vt(\log \tilde{t}-i\pi I)|\Phi_{s}\rangle+
\right.\]
\[+\left. \;\langle\Omega_{r}| t^{\frac{\pi}{\hb}}
\vt(\log t-i\pi I)^*|\Psi_t\rangle
\langle\Psi_t|\Phi_{\tilde{{s}}}\rangle
V_1(\tilde{{s}})|\langle\Phi_{ \tilde{{s}}}|V_1({S_+})
|\Psi_{\tilde{t}}\rangle\langle\Psi_{\tilde{t}}|
\tilde{t}^{\frac{\pi}{\hb}}\;\vt(\log \tilde{t}-i\pi )|\Phi_{s}\rangle
\right\}=\]
\[=\langle\Omega_{r}|\Phi_{s}\rangle^{-1}\;\times\;\left\{\; V_1(-\tilde{{s}})\langle\Omega_{r}|\;
\vt(\log T_+-i\pi I)^*|\Phi_{\tilde{{s}}}\rangle
\langle\Phi_{ \tilde{{s}}}|
\vt(\log T_+-i\pi I)|\Phi_{s}\rangle+\right.\]
\[+\left. \;V_1(\tilde{{s}})\langle\Omega_{r}| T_+^{\frac{\pi}{\hb}}
\vt(\log T_+-i\pi I)^*|\Phi_{\tilde{{s}}}\rangle
\langle\Phi_{ \tilde{{s}}}
| T_+^{\frac{\pi}{\hb}}\;\vt(\log T_+-i\pi I)|\Phi_{s}\rangle
\right\}\ .\]
Therefore by  (\ref{efw}) and  (\ref{efwt}) we obtain
\[V_1({r})V_1(-{s})=i(-1)^r c_\hbar^\prime  
e^{-i\frac{\log^2 \tilde{{s}}}{\hbar}}\langle\Omega_{r}|\Phi_{s}\rangle^{-1}\;\times\]
\[\times\;\left\{\;V_1(-\tilde{{s}})
\langle\Phi_{\tilde{{s}}}|
T_+^\frac{\pi}{\hb}\vt(\log T_+-i\pi I)|\Phi_{r}\rangle
\langle\Phi_{ \tilde{{s}}}|
\vt(\log T_+-i\pi I)|\Phi_{s}\rangle+\right.\]
\beq
\!\!\!\!\!\!+\left. \;V_1(\tilde{{s}})
\langle\Phi_{\tilde{{s}}}|
\vt(\log T_+-i\pi I)|\Phi_{r}\rangle
\langle\Phi_{ \tilde{{s}}}
| t^{\frac{\pi}{\hb}}\;\vt(\log T_+-i\pi I)|\Phi_{s}\rangle
\right\}\label{l1+1-}
\ .\eeq	
On the other hand
\[\langle\Omega_{s}|V_1(-R_+)V_1({S_+})|\Phi_{r}\rangle=\langle\Omega_{s}|X(3,3)|\Phi_{r}\rangle=\]
\[=
\langle\Omega_{s}| T_+^{\frac{\pi}{\hb}}\;\vt(\log T_+-i\pi I)^*\;V_1(-{s})\;
 T_+^{\frac{\pi}{\hb}}\vt(\log T_+-i\pi I)\;|\Phi_{r}\rangle+\]
 \[+\;\langle\Omega_{s}|
\;\vt(\log T_+-i\pi I)^*\;V_1({s})\;
\;\vt(\log T_+-i\pi I) |\Phi_{r}\rangle \ .\]
Hence
\[V_1(-{s})V_1({r})=\langle\Omega_{s}|\Phi_{r}\rangle^{-1}\;\times\]
\[\times\;\left\{\;
\langle\Omega_{s}| t^{\frac{\pi}{\hb}}
\vt(\log t-i\pi )^*|\Psi_t\rangle
\langle\Psi_t|\Phi_{\tilde{{s}}}\rangle
\langle\Phi_{ \tilde{{s}}}|V_1(-{s})
|\Psi_{\tilde{t}}\rangle\langle\Psi_{\tilde{t}}|
 \tilde{t}^{\frac{\pi}{\hb}}\;\vt(\log \tilde{t}-i\pi )|\Phi_{r}\rangle
+\right.\]
\[+\left. \;
\langle\Omega_{s}|\;
\vt(\log t-i\pi I)^*|\Psi_t\rangle
\langle\Psi_t|\Phi_{\tilde{{s}}}\rangle
\langle\Phi_{ \tilde{{s}}}|V_1({s})
|\Psi_{\tilde{t}}\rangle\langle\Psi_{\tilde{t}}|\vt(\log \tilde{t}-i\pi )|\Phi_{r}\rangle
\right\}=\]
\[=\langle\Omega_{s}|\Phi_{r}\rangle^{-1}\;\times\;\left\{\;
V_1(-\tilde{{s}})
\langle\Omega_{s}| T_+^{\frac{\pi}{\hb}}
\vt(\log T_+-i\pi I)^*|\Phi_{\tilde{{s}}}\rangle
\langle\Phi_{ \tilde{{s}}}
| T_+^{\frac{\pi}{\hb}}\;\vt(\log T_+-i\pi I)|\Phi_{r}\rangle
+\right.\]
\[+\left. \;
V_1(\tilde{{s}})\langle\Omega_{s}|\;
\vt(\log T_+-i\pi I)^*|\Phi_{\tilde{{s}}}\rangle
\langle\Phi_{ \tilde{{s}}}|
\vt(\log T_+-i\pi I)|\Phi_{r}\rangle
\right\}\ .\]
Therefore by  (\ref{efw}) and  (\ref{efwt}) 
\[V_1(-{s})V_1({r})=i(-1)^k c_\hbar^\prime  
e^{-i\frac{\log^2 \tilde{{s}}}{\hbar}}\langle\Omega_{s}|\Phi_{r}\rangle^{-1}\;\times\]
\[
\times\;\left\{\;
V_1(\tilde{{s}})
\langle\Phi_{\tilde{{s}}}|
\vt(\log T_+-i\pi I)|\Phi_{r}\rangle
\langle\Phi_{ \tilde{{s}}}
| T_+^{\frac{\pi}{\hb}}\;\vt(\log T_+-i\pi I)|\Phi_{s}\rangle
+\right.\]
\beq
\!\!\!\!\!\!+\left. \;
V_1(-\tilde{{s}})
\langle\Phi_{\tilde{{s}}}|
T_+^\frac{\pi}{\hb}\vt(\log T_+-i\pi I)|\Phi_{r}\rangle
\langle\Phi_{ \tilde{{s}}}|
\vt(\log T_+-i\pi I)|\Phi_{s}\rangle
\right\}
\label{p1+1-}\ .
\eeq
Applying  (\ref{omfi}) we see that (\ref{l1+1-}) and 
 (\ref{p1+1-}) are the same. Thus  
\[V_1({r})V_1(-{s})= V_1(-{s})V_1({r})\]
so (\ref{1+1-}) holds.

The proofs of formulae 
    \rf{1+2-},
 \rf{1-2+} and \rf{2-1+} and \rf{1+2+}
 are exactly the same so we omit them.

In order to prove  \rf{2--} one has to observe additionally that
\[X(1,4)=X(4,1)=0\ .\]
Thus we proved formulae \rf{1++} $\div$
\rf{1+2+}.

Adding  \rf{1++} and \rf{1-2+}
 and substracting \rf{1++} from \rf{1-2+}
 yields
\[V_1(r)\{V_1(s)+V_2(s)\}=V_1(r)V(s,1)=\{V_1(s)+V_2(s)\}V_1(r)=V(s,1)V_1(r).\]
Hence
\beq
\label{p1}
V_1(r)V(s,1)=V(s,1)V_1(r)\mi V_1(r)V(s,-1)=V(s,-1)V_1(r)
\eeq
In the same manner we can see that
\[
V(r,\varrho)V(s,\sigma)=V(s,\sigma)V(r,\varrho)
\]
for any $r,s\in\R\setminus\{0\}$ and $\varrho,\sigma\in\{0,1\}$.
Moreover \rf{2--} holds
\[V_2(r)V_2(-s)=0=V_2(-s)V_2(r)\ .\]\mqed
We stress  that satisfying \rf{wz:przemN} and \rf{2--} 
is necessary, but not sufficient condition  
 for  $V$ to be a representation of  $N$.
%\bwn
%Aby zna/c posta/c wszystkich unitarnych reprezentacji 
%kategorii N w spacei ${\cal K}$ wystarczy zna/c posta/c reprezentacji 
%w $\C$.
%\ewn

We find now a formula for all unitary representations of  $N$ 
acting on $\C$. 
\bst
\label{repNwC} All unitary representations of
$N$ acting on $\C$ are of the form
\[ V(R,\rho)=\fh(M R, \mu \rho
\chi(M R<0))\]
where $M\in \R$ and $\mu =\pm 1$.
\est
\dow
We first find solutions of equation  \rf{2--}.
There are 3 cases
\bit
\item[1.] for any  $r\in\R_+$ we have $\od V_2(r)=0$ and $V_2(-r)\neq 0$.
 
Then:
\[V(r,\varrho)=\left\{\ba{ccc}
V_1(r)& \dla & r\rangle0\\
V_1(r)+\varrho V_2(r)&\dla & r<0
\ea\right.,\]
 so
\[V(R,\rho)=V_1(R)+\rho
\chi(R<0) V_2(R)=V(R,\rho\chi(R<0))=V(\phi(R,\rho)).\]
Similarly
\[V(S,\sigma   )=V_1(S)+\sigma\chi(S<0) 
V_2(S)=V(S,\sigma\chi(S<0))=V(\phi(S,\sigma)).\]
Moreover, since  $V$ is representation of  $N$ 
\[V(\phi(R,\rho))V(\phi(S,\sigma))=
V(\phi((R,\rho)\mlot_N(S,\sigma))).\]
%na mocy Stwierdzenia \ref{an3}. 
%
Hence by Theorem \ref{twstare}
 $$V(R,\rho)=\fh(MR,\mu\rho\chi(R<0)),$$
where $M\geq 0$ and $\mu=\pm 1 $. 

\item[2.] For any $r\in\R_+$ we have  $\od V_2(-r)=0$ and $ V_2(r)\neq 0$.\\
In the same manner as in the previous case we conclude that 
 $$V(R,\rho)=\fh(-MR,\mu\rho\chi(R>0))=
\fh(\tilde{M}R,\mu\rho\chi(\tilde{M}R<0)),$$
where $\tilde{M}=-M\leq 0$ and $\mu=\pm 1 $.

\item[3.] For any $r\in\R_+$ we have  $\od V_2(r)=0\od$ and $\od V_2(-r)=0$.\\ 
Then 
 $V(R,\rho)=V_1(R)$ and $V(S,\sigma)=
V_1(S)$.
Note that
 $$V_1(R)V_1(S)\neq V_1([R+S]
_{\mu\chi(\eh S^{-1}R<0)}),$$
 since $[R+S]_{\mu\chi(\eh S^{-1}R<0)}$ 
depends on $\rho$ and on $\sigma$, while the left hand side does not.
 Therefore we conlude that $V=V_1$ is not a  representation of  $N$. 
\eit
Thus we proved that all unitary representations of  $N$ acting 
 on Hilbert space $\C$ 
are 
\[V(R,\rho)=\fh(MR,\mu \rho\chi(MR<0) ),\]
where $M\in \R$ and $\mu =\pm 1$.\hfill\qed

Now we turn to the case of representations of $N$ acting on arbitrary Hilbert 
space $\cK$.
 Note that  
if $\dim\cK=k<\infty$, then from commutation of  unitary operators 
$V(r,\rho)$ and $V(s,\sigma)$ follows the existence of an 
 ortonormal basis diagonalizing matrices  $V(r,\rho)$ and $V(s,\sigma)$
for any $r,s\in\R$ and $\rho,\sigma\in\{-1,1\}$.  
Thus  the problem reduces to finding solutions of  $k$ scalar equations 
\[V_o(R,\rho)V_o(S,\sigma)=V_o((R,\rho)\mlot_N (S,\sigma))\ ,\]
where complex-valued function $V_o$ is defined on  $\R\times \{-1,1\}$.
The same conclusion can be drawn for an arbitrary separable 
 Hilbert space $\cK$.  
 The reason is that  operators $V(r,\rho)$ and $V(s,\sigma)$ 
belong to commutative *-subalgebra  of  $B(\cK)$. 
Therefore, by spectral  theorem 
 and its consequences \cite[Chapter X]{dunfII} 
 operators $V(r,\rho)$ and $V(s,\sigma)$ have the same spectral measure
\[V(r,\rho)=\int_0^{2\pi}V_o(r,\rho,t)dE_{\cK}(t)\mi 
V(s,\sigma)=\int_0^{2\pi}V_o(s,\sigma,t)dE_{\cK}(t).\]
Hence 
\[V(R,\rho)=\int_{\R\times \{-1,1\}} 
\int_0^{2\pi}V_o(r,\rho,t)dE_{\cK}(t)\te dE_{R,\rho}(z)\]
and
\[V(S,\sigma)=\int_{\R\times \{-1,1\}} \int_0^{2\pi}V_o(s,\sigma,t)
dE_{\cK}(t)\te dE_{S,\sigma}(z)\ ,\]
where $dE_{R,\rho}$ is the joint spectral measure of strongly commuting 
operators $R$ and $\rho$.

Thus theorem 
 \ref{twN} reduces to 
 the already proved  Proposition 
\ref{repNwC}.
\hfill\qed

\section{The braided quantum group  $M$ } 
\label{mul}
The main goal of this paper is finding all unitary representations of the operator domain $M$, which will be introduced below.
 This operator domain is very close to the operator domain correponding  
to the quantum  'ax+b' group.  To emphasize that 
 we use the same letters  $b$ and $\beta$,  which denoted operators  
generating \footnote{In fact,  generators are $b$ and $ib\beta$} 
quantum 'ax+b' group in \cite{ax+b}.  
 Theorem \ref{twM} we prove in this section will be crucial  
 in our next paper \cite{paper2}.

Consider operators $b$ 
 and  $\beta$  such that
\beq
\label{rkmul}
b=b^*\mand \beta=\beta^*\mand \beta b=-b\beta\mand 
\beta^2=\chi(b\neq 0) \ .\eeq
Define an operator domain $M$ by
\[M_{\cH}=\{(b,\be )\in\ch ^2|b=b^*,\be =\be ^*,
\be  b=b\be ,\be ^2=\chi(b\neq 0)\}\ \]
 and $M^2$ by  
\[M_{\cH}^2=
\{\left((b,\be),(d,\de)\right)|
(b,\be ),(d,\delta)\in M_{\cH}\;\;b\za d, \;\; b\delta=\delta b,  
\;\; d\beta=\beta d,\;\;
\beta\delta=\delta\beta \}\ .\]
%
%Note that $M^2$ is also an operator domain.

The  colorally below allows us to define group operation on $M$.
\bfa 
Let  $\left( (b,\beta),(d,\delta)\right)\in M_{\cH}^2 $ and 
let $\ker d=\{0\}$. Let
\[f=e^{\frac{i\hbar}{2}}d^{-1}b\mand\phi=\pm \beta\delta \chi(\eh d^{-1} b<0)\]
and
\[
[b+d]_{\phi }=F_{\hbar}(f,\phi)^* d
F_{\hbar}(f,\phi)\mand
\tilde{\delta}=F_{\hbar}(f,\phi)^*\delta 
F_{\hbar}(f,\phi)\ .\]
Then $$([b+d]_{\phi},\tilde{\delta})\in M_{\cH}.$$
\efa
\dow We first show that  
$[b+d]_{\phi}
%=F_{\hbar}(f,\phi)^* dF_{\hbar}(f,\phi)
$  is a selfadjoint operator. Operator $\phi$ is a 
reflection operator corresponding to  
$e^{\frac{i\hbar}{2}}d^{-1}b$,
 since $\phi$  is selfadjoint, $\phi^2=\chi(e^{\frac{i\hbar}{2}}d^{-1}b<0)$ 
and 
$\phi$ anticommutes with $b+d$. Operator $b+d$ is selfadjoint, if 
$\eh b d \geq 0$ (see Theorem 5.4  \cite{qexp}).
Hence $[b+d]_{\phi}=F_{\hbar}(f,\phi)^* d
F_{\hbar}(f,\phi)$  is an  selfadjoint extension of 
sum of  selfadjoint operators  $\eh f d$ and $d$ (see Theorem 6.1 and 
4.1 \cite{qexp}), corresponding to the reflection  operator $\phi$.

We prove now that  $\tilde{\delta}$ is a  selfadjoint operator.
From
\[\tilde{\delta}=F_{\hbar}(f,\phi)^*\delta F_{\hbar}(f,\phi)\]
follows that $\tilde{\delta}$ is selfadjoint, since it  is unitarily 
 equivalent to selfadjoint operator $\delta$.
Moreover, note that    $[b+d]_{\phi}$  commutes with $\tilde{\delta}$, 
because $d$ commutes with  $\delta$. 

In order to prove that $([b+d]_{\phi},\tilde{\delta})\in M_{\cH}$,  
 we have to check that
\[\tilde{\delta}^2=\chi ([b+d]_{\phi}\neq 0)\ .\]
To this end compute 
 \[\tilde{\delta}^2=F_{\hbar}(f,\phi)^*\delta 
F_{\hbar}(f,\phi)F_{\hbar}(f,\phi)^*\delta 
F_{\hbar}(f,\phi)=
F_{\hbar}(f,\phi)^*\delta^2 
F_{\hbar}(f,\phi)=\]
\[=F_{\hbar}(f,\phi)^*\chi(d\neq 0)F_{\hbar}(f,\phi)=
\chi(F_{\hbar}(f,\phi)^*d F_{\hbar}(f,\phi)\neq 0)=\chi([b+d]_{\phi}\neq 0)\]    \mqed 
We are now in a position to define operation on  $M$ by
\[\mlot_M:M^2_{\cH}\longrightarrow M_{\cH}\]
\[(b,\be )\mlot_M (d,\delta)=([b+d]_{\phi},\tilde{\delta}).\]
Operation $\mlot_M$ is associative. Moreover, $M$ with this operation 
forms a braided quantum group.
\subsection{The matrix representation of $M$}
Let $\ker b=\ker d=\{0\}$. 
Since  $b\;\za\; d$
 it follows that  $b$ and $d$ commute with $\sign\; b$ and $\sign\; d$. 
Therefore we may introduce notation
\[\ba{lll}
 {\cH}_{++}&=&{\cH}(b>0)\cap {\cH}(d>0)\\
 {\cH}_{+-}&=&{\cH}(b>0)\cap {\cH}(d<0)\\ 
 {\cH}_{-+}&=&{\cH}(b<0)\cap {\cH}(d>0)\\
 {\cH}_{--}&=&{\cH}(b<0)\cap {\cH}(d<0)\\
\ea\ .\]
Then  ${\cH}={\cH}_{++}\oplus {\cH}_{+-}\oplus {\cH}_{-+}\oplus {\cH}_{--}$.
Every vector  $\psi\in{\cH}$ we represent by 
\[\psi=
\left[
\ba{ccc}
\psi_{++}\\
\psi_{+-}\\
\psi_{-+}\\
\psi_{--}
\ea\right],\]
where $\psi_{++}\in {\cH}_{++}$, $\psi_{+-}\in {\cH}_{+-}$, 
$\psi_{-+}\in {\cH}_{-+} $ and $\psi_{--}\in {\cH}_{--} $.
Therefore operators acting on ${\cH}$ will be represented 
by  $4\times 4$ matrices.
From (\ref{rkmul}) follows that maps 
$\be : {\cH}_{++} \rightarrow {\cH}_{-+}\mand \be : 
{\cH}_{-+} \rightarrow {\cH}_{++} $
and 
$\be  : {\cH}_{--}\rightarrow {\cH}_{+-} \mand  \be  : 
{\cH}_{+-}\rightarrow {\cH}_{--}$
are mutually inverse.
Similarly by (\ref{rkmul}) we obtain   that maps 
$\de : {\cH}_{++}\rightarrow {\cH}_{+-}\mand \de :
{\cH}_{+-}\rightarrow {\cH}_{++}$
and
$\de  : {\cH}_{--}\rightarrow {\cH}_{-+} \mand  \de  : 
{\cH}_{-+}\rightarrow {\cH}_{--}$
are mutually inverse.
Hence Hilbert spaces 
${\cH}_{++},{\cH}_{-+},{\cH}_{+-}$ and ${\cH}_{--}$ 
are unitarily equivalent. 
In what follows for simplicity we assume  
that ${\cH}_{++}={\cH}_{-+}={\cH}_{+-}={\cH}_{--}$ and denote  
this Hilbert space by ${\cH}_+$. 
With this notation we have the  
following  representations of  $\beta$ and $\delta$
\[
\beta=\left[\begin{array}{cccc}
0&0&I&0\\0&0&0&I\\I&0&0&0\\0&I&0&0
\end{array}\right],\hspace{20mm}\mand 
\delta=\left[\begin{array}{cccc}
0&I&0&0\\I&0&0&0\\0&0&0&I\\0&0&I&0
\end{array}\right].
\]
Since $b$ anticommutes with  $\beta$ and commutes with $\delta$ 
and  $d$  anticommutes with $\de  $ and commutes with $\be   $, 
they will be represented as follows
\[b=\left[\begin{array}{cccc}
b_+&0&0&0\\0&b_+&0&0\\0&0&-b_+&0\\0&0&0&\hspace{-2mm}-b_+
\end{array}\right],\hspace{10mm}\mand
d=\left[\begin{array}{cccc}
d_+&0&0&0\\0&-d_+&0&0\\0&0&d_+&0\\0&0&0&-d_+
\end{array}\right],
\]
where   $b_+$ and $d_+$ are  restrictions of  $b$ and $d$ to ${\cH}_+$. 
It is easily seen that   $b_+$ and $d_+$ are selfadjoint and 
strictly positive and    $b_+\za \;d_+$.
Moreover
\[
f=\eh d^{-1}b=\left[\begin{array}{cccc}
f_+&0&0&0\\0&-f_+&0&0\\0&0&-f_+&0\\0&0&0&f_+
\end{array}\right],
\]
where   $f_+=\eh d_+^{-1}b_+$
 is selfadjoint and 
strictly positive.
Hence
\[
\phi= (-1)^k\beta\delta \chi(\eh d^{-1}b<0)=
 (-1)^k  \beta\delta \chi(f<0)=
\left[\begin{array}{cccc}
0&0&0&0\\0&0&(-1)^k I&0\\0&(-1)^k I&0&0\\0&0&0&0
\end{array}\right],
\]
where $k\in \N$.
\subsection{From $N$ to $M$}
As we said in Introduction, the operator domain $N$ is auxiliary, what we 
 are really interested in are unitary representations of $M$.  
However,  $N$ was easier to work with, because it was commutative,
  so we found formula for all 
unitary representations of it.
 It would be very nice now to have an  operator map from $N$ into $M$, which would allow us to ``tranfer'' results obtained for $N$ to $M$.
Such an operator map is constructed below
\bst
\label{n>m}
%Let  ${\cal H}$ b/edzie  spacei/a Hilberta.
Let  
   $((R,\rho),(S,\sigma))\in N  _{\cal H}^2$.
Define a map
$$\fil:N_{\cal H}\rightarrow M_{\C^2\te{\cal H}}$$ 
for any $(R,\rho )\in N_{\cal H}$  by  
\[\fil(R,\rho )=(\left[\ba{cc} R&0\\0&-R\ea\right],
\left[\ba{cc}0&  \rho\\\rho&0\ea\right]).\]
Then $\phi$ is an  operator map and 
%\bit
%\item[1.]
%{$ \fil (R,\rho )\in M_{\C^2\te {\cal H}}$ }
%\item[2.]
%{Odwzorowanie 
%\[\fil: N_{\cal H}\longrightarrow M_{\C^2\te {\cal H}}\]
%komutuje z do/l/aczonym dzia/laniem izometrii.}
%\item[3.]
 for any $\left((R,\rho ),(S,\sigma)\right)\in N_{\cal H}^2 $ we have
%\[\left(\fil (R,\rho ),\fil
%(S,\sigma)\right)\in M_{\C^2\te {\cal H}}^2  \]}
%\item[4.]
\[\fil  \left((R,\rho )\mlot_{N_{\cal H}} (S,\sigma)\right)=
 \fil (R,\rho ) \mlot_{M_{\C^2\te {\cal H}}} \fil
(S,\sigma).\] 
%}\eit
\est
%%%%% 
%W dowodzie pominiemy trywialne punkty 1,2,3.
%Udowodnimy teraz punkt 4.
%$ \left[\ba{cc} R&0\\0&-R\ea\right]$ and 
%$\left[\ba{cc}0&  \rho\\\rho&0\ea\right]       $ are samosprz/e/zone i
% antykomutuj/a ze sob/a  and 
%\[\left[\ba{cc}0&  \rho\\\rho&0\ea\right]^2=
%\left[\ba{cc} I&0\\0&I\ea\right]
%=\chi \left(\left[\ba{cc} R&0\\0&-R\ea\right]\neq 0\right)\]
%czyli $\fil (R,\rho )\in M_{\C^2\te {\cal H}}$ }
%\item[2.]{ Bezpo/srednim rachunkiem mo/zna sprawdzi/c, /ze dla 
%dowolnego $\hat{V}\in {\rm Unit}({\cal H})$:
%\[\left[\ba{cc} \hat{V}^*&0\\0&\hat{V}^*\ea\right]\left(\left[\ba{cc}R&0\\0&R
%\ea
%\right],\left[\ba{cc}0&\rho\\\rho&0\ea\right]\right)
%\left[\ba{cc} \hat{V}&0\\0&\hat{V}\ea\right]=\]
%\[=
%\left(\left[\ba{cc}\hat{V}^*R\hat{V}&0\\0&\hat{V}^*R\hat{V}
%\ea
%\right],\left[\ba{cc}0&\hat{V}^*\rho\hat{V}\\
%\hat{V}^*\rho\hat{V}&0\ea\right]\right),\]
%czyli  $\fil$ komutuje z 
%do/l/aczonym dzia/laniem izometrii.}
%\item[3.]{/Latwo sprawdzi/c, /ze spe/lnione are wszystkie warunki zgodno/sci:
%\[\left[\ba{cc} R&0\\0&-R\ea\right]\za\;\left[\ba{cc} S&0\\0&-S\ea\right] \]
%i
%
%\[ \left[\ba{cc} R&0\\0&-R\ea\right] \left[\ba{cc}0&  \sigma\\\sigma&0\ea\right]      =\left[\ba{cc}0&  -\sigma R\\\sigma R&0\ea\right]=\left[\ba{cc}0&  \sigma\\\sigma&0\ea\right]\left[\ba{cc} R&0\\0&-R\ea\right]\] and 
%\[ \left[\ba{cc} S&0\\0&-S\ea\right] \left[\ba{cc}0&  \rho\\\rho&0\ea\right]
 %     =\left[\ba{cc}0&  -\rho S\\\rho S&0\ea\right]=\left[\ba{cc}0&  \rho\\\rho&0
%\ea\right]\left[\ba{cc} S&0\\0&-S\ea\right]\]}
%\item[4.]{Jedyna nietrywialna rzecz w tym dowodzie. 
\dow We first prove that
 \[\left[\left[\ba{cc} R&0\\0&-R\ea\right]+\left[\ba{cc} S&0\\0&-S\ea\right]
\right]_{\phi}  =\left[\ba{cc} [R+S]_{\tau}&0\\
0&-[R+S]_{\tau}\ea\right].\]
Since selfadjoint extensions are determined uniquely by 
 reflection operators, it is enough to check that $I_{\C^2}\te \tau = \phi $.
We already know that
\[\phi= (-1)^k \beta\delta=(-1)^k \left[\ba{cc}0&  \rho\\\rho&0\ea\right]
 \left[\ba{cc}0&  \sigma\\\sigma&0\ea\right]= 
(-1)^k \left[\ba{cc}  \rho\sigma&0\\0&\rho\sigma\ea\right]\]
 and
 \[\tau =(-1)^k
\rho\sigma \ .\]
Hence
%, dla ustalonego  $k\in \N$,
\[\phi  =\left[\ba{cc}  \tau&0\\0&\tau\ea\right].\]
Therefore  selfadjoint extensions given by this  
  reflection operators 
are the same.

It remains to prove that
\[\fil(\tilde{\sigma})=\tilde{\delta}\ .\] 
To this end, note that
\[
\fil(\tilde{\sigma})=
\left[
\ba{ccc}
&0 & \tilde{\sigma} \\ 
&\tilde{\sigma}& 0 \\
\ea
\right]
=
\]
\[=\left[\ba{cc}0&
\fh(T,\tau\chi(T <0))^*\sigma 
\fh(T ,\tau\chi(T <0))\\ 
\fh(T ,\tau\chi(T <0))^*\sigma 
\fh(T ,\tau\chi(T <0))&0\ea\right]=\]
\[=\fh(f,\phi)^*\delta 
\fh(f,\phi)=\tilde{\delta}\ .\]\mqed
\subsection{Untary representations of   $M$}
We can now formulate our main result.
\btw
 \label{twM}
$U$ is an unitary  representation of  $M$ 
acting on Hilbert space ${\cal K}$, if there 
 exists such $(g,\gamma)\in M_{\cal K}$  that for any $(b,\beta)
% ,d,\delta)
\in M_{\cal H}$ 
\beq
\label{repM}
U( b,\beta)=\fh(g\te b,(\gamma\te \beta)
\chi(g\te b<0)) .
\eeq
\etw
\dow We first show that if $U$ is a representation of 
$M$, it has form \rf{repM}. 
Let $\fil:N_{\cal H}\rightarrow M_{\C^2\te{\cal H}}$ 
be the  operator map considered before.
Let  $(b,\beta)$ denote the following  element from  
 $M_{\C^2\te{\cal H}}$
\[
(b,\beta):=\fil({R},\rho)=\left(\left[\ba{cc}R&0\\0&-R\ea\right], 
\left[\ba{cc}0&\rho\\ \rho&0\ea\right]\right)
\ .\]
Let
\[V(R,\rho):=U(b,\beta)=U(\fil(R,\rho))\ .\]
Since  $U$ is a unitary representation of  
$M$, it follows by Proposition  \ref{n>m} that  
$V$ is a unitary representation of  
$N$ on ${\cal K}\te\C^2 $.
Next by  Theorem \ref{twN} we get
\beq
\label{repModK}
V(R,\rho)=U\left(\left[\ba{cc}R&0\\0&-R\ea\right], 
\left[\ba{cc}0&\rho\\ \rho&0\ea\right]\right)=
\fh(M\te R,(\mu\te\rho)\chi(M\te R<0))
\eeq
where $(M,\mu)\in N_{{\cal K}\te\C^2 }$.
 We find the conditions on  
$(M,\mu)\in N_{\C^2\te {\cal K} }$, under which  \rf{repModK} holds.
We know that $U,\fil$ and $V$ are  operator maps. 
Let a unitary operator 
${\hat{{\cal U}}}$ 
 be given by
 \[\hat{{\cal U}}:=I_{\cal K}\te \left[\ba{cc}0&1\\
1&0\ea\right]\te I_{\cal H}\ .\]
Apply  ${\rm ad}_{\hat{{\cal U}}}$ to the  left hand side of \rf{repModK}. 
We obtain
 \[\left( I_{\cal K}\te \left[\ba{cc}0&1\\
1&0\ea\right]\te I_{\cal H}  \right)^*U\left(\left[\ba{cc}R&0\\0&-R\ea\right], 
\left[\ba{cc}0&\rho\\ \rho&0\ea\right]\right)
\left( I_{\cal K}\te \left[\ba{cc}0&1\\
1&0\ea\right]\te I_{\cal H}  \right)
=\]
\[=U\left(\left[\ba{cc}-R&0\\0&R\ea\right], 
\left[\ba{cc}0&\rho\\ \rho&0\ea\right]\right)\ .\]
Hence 
\beq
\label{rozdzielanie}
\hat{{\cal U}}^*\fh(M\te R,\mu\te\rho)\hat{{\cal U}}=
\fh(-M\te R,\mu\te\rho)\ .
\eeq
Observe that for any $t\in\R\setminus\{0\}$ the pair $(tR,\rho)$ 
belongs to $\N_{\cH}$ if only $(R,\rho)$ does. Therefore 
we may put $tR$ instead of $R$ in \rf{rozdzielanie}.
The function $\fh$ is not injective, however the family 
$\fh (t\cdot)$ separates  points of $\Delta_{\rm real}$. 
Therefore from \rf{rozdzielanie} follows
\[\left(I_{\cal K}\te \left[\ba{cc}0&1\\
1&0\ea\right]\right)M\left(I_{\cal K}\te \left[\ba{cc}0&1\\
1&0\ea\right]\right)=-M\]
and
\[\left(I_{\cal K}\te \left[\ba{cc}0&1\\
1&0\ea\right]\right)\mu\left(I_{\cal K}\te \left[\ba{cc}0&1\\
1&0\ea\right]\right)=\mu .\]
Let us introduce another unitary operator
$$\hat{{\cal V}}:=I_{\cal K}\te \left[\ba{cc}1&0\\
0&-1\ea\right]\te I_{\cal H}\ .$$
Applying ${\rm ad}_{\hat{{\cal V}}}$
 to the left hand side of \rf{repModK}
we obtain
 \[\left( I_{\cal K}\te \left[\ba{cc}1&0\\
0&-1\ea\right]\te I_{\cal H}  \right)^*
U\left(\left[\ba{cc}R&0\\0&-R\ea\right], 
\left[\ba{cc}0&\rho\\ \rho&0\ea\right]\right)
\left( I_{\cal K}\te \left[\ba{cc}1&0\\
0&-1\ea\right]\te I_{\cal H}  \right)
=\]
\[=U\left(\left[\ba{cc}R&0\\0&-R\ea\right], 
\left[\ba{cc}0&-\rho\\ -\rho&0\ea\right]\right)\ .\]
Hence 
\[
\hat{{\cal V}}^*\fh(M\te R,\mu\te\rho)\hat{{\cal V}}=
\fh(M\te R,-\mu\te\rho),
\]
so
\[\left(I_{\cal K}\te \left[\ba{cc}1&0\\
0&-1\ea\right]\right)M\left(I_{\cal K}\te \left[\ba{cc}1&0\\
0&-1\ea\right]\right)=M\]
and
\[\left(I_{\cal K}\te \left[\ba{cc}1&0\\
0&-1\ea\right]\right)\mu\left(I_{\cal K}\te \left[\ba{cc}1&0\\
0&-1\ea\right]\right)=-\mu .\]

Therefore $M$ and $\mu$ have form
\[M=g\te \left[\ba{cc}1&0\\
0&-1\ea\right]\mi \mu=\gamma\te \left[\ba{cc}0&1\\
1&0\ea\right],\]
where $g$ and $\gamma$ are operators acting on Hilbert space ${\cal K}$.
Note that since $(M,\mu)\in N_{{\cal K}\te \C^2}$, we see that  
$g,\gamma$ are selfadjoint, $g$ anticommutes with $\gamma$ and $\gamma^2 =
\chi(g\neq 0)$. It means that $(g,\gamma)\in M_{\cal K}$.
Hence
\[
U\left(\left[\ba{cc}R&0\\0&-R\ea\right], 
\left[\ba{cc}0&\rho\\ \rho&0\ea\right]\right)=\]
\[=\fh(g\te\left[\ba{cc}R&0\\0&-R\ea\right] ,
(\gamma\te\left[\ba{cc}0&\rho\\ \rho&0\ea\right])
\chi(g\te\left[\ba{cc}R&0\\0&-R\ea\right] <0))
,\]
so
\[
U(b,\beta)=\fh(g\te b,(\gamma\te\beta)
\chi(g\te b<0))
,\]
where $(g,\gamma)\in M_{\cal K}$.

We proved that every unitary representation of $M$ has form \rf{repM}.
What is left is to show that every operator  function given by \rf{repM}
 is  a unitary representation of $M$.
To this end it is sufficient to show 
\bst
\label{m>a}
Let $(g,\gamma)\in M_{\cal K}$ and 
 $\left((b,\beta),(d,\delta)\right)\in M_{\cal H}^2$.
For any $(b,\beta)\in M_{\cal H}$ define a map by  
\[\fil_{(g,\gamma)}: M_{\cal H}\ni(b,\beta)\mapsto 
(g\te b, \gamma\te \beta)\in N_{{\cal K}\te {\cal H}}\]
Then
\[\fil_{(g,\gamma)} \left((b,\beta)\mlot_{M} (d,\delta)\right)=
 \fil_{(g,\gamma)}(b,\beta) \mlot_{N} \fil_{(g,\gamma)}
(d,\delta)\ .\] 
\est
\dow  We first have to check that selfadjoint extensions on both 
sides of the above formula are the same, i.e. 
that  $g\te [b+d]_{\phi}$  equals $[g\te b + g\te d]_{\tau}$.
Since  $ [g\te b + g\te d]_{\tau}=g\te [b +  d]_{{\tau |_{\cal H}}}$ 
and selfadjoint extensions are given uniquely by 
reflection operators, it is enough  to check that $\phi = \tau |_{\cal H}$.
We know that
\[\phi= (-1)^k \beta\delta
\chi(\eh d^{-1}b <0)\]
 and
 \[\tau =(-1)^k
(\gamma\te \beta)(\gamma\te \delta)\chi (\eh (g\te b)^{-1}(g\te d)<0)=\]
\[ =(-1)^k
(\gamma^2\te \beta\delta)\chi (I\te\eh b^{-1} d<0)\ .\]
Since $\gamma^2=\chi(g\neq 0)$, it follows that 
\[\tau =(-1)^k
(I\te \beta\delta)\chi (I\te\eh b^{-1} d<0)\ .\]
Obviously
\[\tau|_{\cal H}= (-1)^k \beta\delta
\chi(\eh d b <0)=\phi,\]
so   
 selfadjoint operators determined by $\tau$ and $\phi$ are also the same.

Secondly, we have to prove that
\[\widetilde{(\gamma\te\delta)}
=(\gamma\te\tilde{\delta})\ .\]
Left hand side is by definition
\[LHS =\fh(I_{\cal K}\te f,I_{\cal K}\te\phi)^*
 (\gamma\te\delta) 
\fh(I_{\cal K}\te f,I_{\cal K}\te\phi)=\]
\[=\left(I_{\cal K}\te\fh(f,\phi)\right)^*
 (\gamma\te\delta)  
\left(I_{\cal K}\te\fh(f,\phi)\right)=RHS\ ,\]
which completes the proof of Proposition \ref{m>a}.
\hfill\qed\\
This finishes also the proof of  Theorem \ref{twM}. 
\hfill\qed\\
This result will prove extremely useful in our next paper 
\cite{paper2}, where we find all unitary representations of the 
quantum 'ax+b' group. 
We also use formula \ref{pom2} derived in this paper in \cite{paper4}. 
\section*{Acknowledgments}
This is part of the author's Ph.D. thesis \cite{phd}, written under the supervision of Professor Stanis{\l}aw L. Woronowicz at the Department of Mathematical Methods in 
Physics, Warsaw University. The author is greatly  
indebted to Professor Stanis{\l}aw L. Woronowicz for stimulating discussions and important hints 
 and comments. The author also wishes to thank 
Professor Wies{\l}aw Pusz and Professor Marek  Bozejko for
several  helpful suggestions.

\end{document}